\input ./style/arxiv-ba.cfg
\documentclass[ba,linksfromyear,preprint]{imsart}
\makeatletter
   \@ifpackageloaded{natbib}{}{\usepackage{natbib}}
\makeatother
\usepackage{mathbh}

\pubyear{2015}
\volume{10}
\issue{4}
\firstpage{877}
\lastpage{908}
\doi{10.1214/15-BA966SI}

\newcommand{\G}{\mathcal{G}}
\newcommand{\jdg}[1]{J_{#1}}
\newcommand{\var}[1]{\mathrm{Var}\left[#1\right]}
\newcommand{\M}[1]{\mathcal{M}(#1)}
\newcommand{\Mno}{\mathcal{M}}
\newcommand{\R}{\mathcal{R}}
\newcommand{\E}[1]{\mathrm{E}\left[#1\right]}
\newcommand{\Ef}[2]{\mathrm{E}_{#1}\left[#2\right]}
\newcommand{\cov}[2]{\mathrm{Cov}\left[#1,#2\right]}
\newcommand{\bel}{y}
\newcommand{\dat}{z}
\newcommand{\indicator}{\mathbh{1}}
\newcommand{\class}[1]{C_{#1}}
\newcommand{\enquote}[1]{``#1''}

\newtheorem{thm}{Theorem}

\begin{document}

\begin{frontmatter}
\title{Posterior Belief Assessment: Extracting Meaningful Subjective
Judgements from Bayesian Analyses with Complex Statistical Models}
%
%
%
\runtitle{Posterior Belief Assessment}

\begin{aug}





\author[addr1]{\fnms{Daniel} \snm{Williamson}\corref{}\ead
[label=e1]{d.williamson@exeter.ac.uk}\ead[label=u1,url]{https://emps.exeter.ac.uk/mathematics/staff/dw356}}
\and
\author[addr2]{\fnms{Michael} \snm{Goldstein}\ead
[label=e2]{michael.goldstein@durham.ac.uk}}

\runauthor{D. Williamson and M. Goldstein}

\address[addr1]{Statistical Science
College of Engineering, Mathematics and Physical Sciences, University
of Exeter, Exeter, \printead{e1}, \printead{u1}}

\address[addr2]{Department of Mathematical Sciences
Durham University, Durham, UK,\\ \printead{e2}}

\end{aug}

%
\begin{abstract}
In this paper, we are concerned with attributing meaning to the results
of a Bayesian analysis for a problem which is sufficiently complex that
we are unable to assert a precise correspondence between the expert
probabilistic judgements of the analyst and the particular forms chosen
for the prior specification and the likelihood for the analysis. In
order to do this, we propose performing a finite collection of
additional Bayesian analyses under alternative collections of prior and
likelihood modelling judgements that we may also view as representative
of our prior knowledge and the problem structure, and use these to
compute posterior belief assessments for key quantities of interest. We
show that these assessments are closer to our true underlying beliefs
than the original Bayesian analysis and use the temporal sure
preference principle to establish a probabilistic relationship between
our true posterior judgements, our posterior belief assessment and our
original Bayesian analysis to make this precise. We exploit second
order exchangeability in order to generalise our approach to situations
where there are infinitely many alternative Bayesian analyses we might
consider as informative for our true judgements so that the method
remains tractable even in these cases. We argue that posterior belief
assessment is a tractable and powerful alternative to robust Bayesian
analysis. We describe a methodology for computing posterior belief
assessments in even the most complex of statistical models and
illustrate with an example of calibrating an expensive ocean model in
order to quantify uncertainty about global mean temperature in the real ocean.
\end{abstract}

%


\end{frontmatter}


\section{Introduction}
The idea that uncertainty is a subjective property of individuals
underlies the development of the field of Bayesian statistics \citep
[\xch{cf.}{c.f.}][]{savage77, lindley00}. You are uncertain about some aspect of
the world, you specify this uncertainty in the form of a prior
probability distribution. You specify a probability model describing
the data generating process that you will use to learn about the world
and update your prior probabilities to posterior probabilities, upon
observation of the data, using Bayes theorem. We term this the
subjective Bayesian approach. \par
Many aspects of a Bayesian analysis are challenging, in particular,
when a prior and likelihood have been established, computing or
sampling from the posterior distribution can be extremely difficult.
This challenge, however, has been well met by the Bayesian field, with
modern techniques in Markov Chain Monte Carlo (MCMC) allowing even the
most complicated statistical models to be sampled from in order to
obtain posterior probabilities \citep[see, for
example,][]{liangetal10}. However, given the complexity of the
statistical models that we are able to develop, the demands of the
science that they apply to, and the power of our computational methods
for sampling from them, methodology for the elicitation of subjective
prior probabilities lags behind and remains extremely difficult and
highly controversial.
\par For these reasons and, perhaps, others, the field of so called
``objective'' Bayesian analysis (or O-Bayes), is now extremely popular
and represents, arguably, the most popular approach to Bayesian
inference. O-Bayes seeks to develop and use automatic prior
distributions for any probability model that, in some sense depending
on the type of prior developed, have a minimal influence on the output
of a Bayesian analysis. For details of the approach see \cite
{berger06}, for a discussion of the competing philosophies between
O-Bayes and subjective Bayes see \cite{goldstein06, berger06} and the
following discussion. This paper concerns the development of the
subjective Bayesian approach, hence we do not comment further on
O-Bayes solutions or philosophies in the rest of our narrative. \par
The statistical models we are now able to develop would, if we were to
attempt a fully subjective Bayesian analysis, often require elicitation
of high dimensional joint prior distributions, perhaps over
spatio-temporal fields or over parameters of multiple types, and
perhaps involving non-standard distributional forms. Current
state-of-the-art elicitation frameworks and tools such as SHELF \citep
{oakleyohagan10} and MATCH \citep{morrisetal14}, enable the elicitation
of univariate distributions of standard forms such as Normal, Beta and
Gamma distributions and do not come close to meeting the requirements
made for some statistical models. Even univariate parameters in a
complex hierarchical statistical model can be subject to so many layers
of abstraction that even understanding what they mean or how they might
impact upon an analysis can make elicitation using current tools
challenging or even dubious. Even if elicitation methodology were
substantially more advanced, time and budget constraints in many
investigations may prohibit the sort of elicitation that involves
either the statistician, the experts or both believing every
probability statement made in order to facilitate the Bayesian analysis.
\par
Our argument then is that even if elicitation methodology catches up
with the computational capabilities of the field and the demand of the
problems we study, in most cases we will not hold any prior probability
distribution or even the likelihood as representative of our beliefs in
complex problems. What meaning can we then give to any Bayesian
analysis that we do perform? Most importantly, is there anything we can
conclude about our own uncertainty from a Bayesian analysis?
\par
Bayes linear methods \citep{goldsteinwooff07} are an offered solution
to this problem. The idea is based on making only partial subjective
prior specifications in the form of means, variances and covariances,
and using the geometry of Hilbert spaces to update these beliefs using
data. Whilst attractive and computationally tractable in many
applications, these methods do not allow us to take advantage of the
many advanced technologies and benefits of the modern fully
probabilistic Bayesian approach, nor do they allow us to combine and
update any fully probabilistic judgements for certain aspects of our
statistical model that we might hold. We also comment that the fully
Bayesian approach is well established and it is therefore easier to
explain to and ultimately publish with collaborators from other fields.
\par
Robust Bayes was a potential avenue of study into this problem. Popular
in the 1980s and 1990s, robust Bayes looked to explore whether any
conclusions from a Bayesian analysis were robust to classes of prior
and likelihood choices \citep[see][for a nice overview]{berger94}.
However, the method focussed on analytic solutions in relatively simple
Bayesian models and, as the computational capabilities of the field
increased exponentially and Bayesian models became more complex, robust
Bayesian approaches became intractable and are now rarely pursued.
\par
Robust Bayes can be thought of as a formalism of sensitivity analysis,
whereby alternative Bayesian analyses are considered to explore how
sensitive posterior inferences are to prior modelling choices \citep
[see, for example,][chapter 5]{gelmanetal04}. Though sensitivity
analysis is an important and useful step in any Bayesian analysis,
there is no formal mechanism for arriving at posterior judgement
following it. For example, what do samples from any alternative
posterior distributions obtained during the sensitivity analysis say
about your actual judgements? Do these alternative posteriors, or
indeed that which formed the main analysis, represent subjective
probabilities? In addition to this, there is no formal method for
ensuring a sensitivity analysis is ``complete'' in the sense that it
fully explores all possible alternative analyses that might have been
deemed reasonable by the analyst and expert.
\par
In this paper, we describe a new methodology, which we term \textit
{posterior belief assessment}, that aims to improve a full Bayes
analysis that we have performed (or can perform) with prior
distributions and likelihood (judgements) that represent our best
current judgement (without our necessarily believing every probability
statement made by these judgements). By best current judgement here, we
mean that we have expert probabilistic judgements that we are unable to
adequately express and that the chosen combination of prior and
likelihood represent these in some way but, due to the complexity of
the problem, we are unable to assert a precise correspondence between
the chosen judgements and our actual beliefs.
\par
Our methodology improves our Bayesian analysis through performing
further calculations under alternative judgements, in order to get
closer to our actual posterior beliefs for key quantities of interest
in a measurable way. Similar to robust Bayes, we attempt to consider
all possible alternative forms that we might give to these judgements.
However, our approach is to use this information along with
foundational arguments, to derive our subjective judgements for key
quantities rather than to look for any posterior probability statements
that are robust to these choices.
\par
Before proceeding, we comment that Bayesian model averaging \citep
[BMA,][]{hoetingetal99}, makes similar arguments regarding the
existence of alternative judgements (or models). That approach, which
puts a probability distribution over a class of alternative models, is,
however, fundamentally different to what is described here. By putting
a subjective probability distribution over any alternative judgements
we consider, as in BMA, we would, under a subjectivist definition of
probability, be making the extremely strong statement that there is a
true model within the set considered. We already know that this is
false in problems of reasonable complexity, such as that given above.
Our approach will not make assumptions of this type. Instead, we begin
with the foundations of subjective inference to establish a
relationship between our underlying judgements and our Bayesian
analysis (Section \ref{meaning}). We then use this relationship to
develop a methodology for combining an original analysis with
alternatives to get closer to those judgements (Sections \ref
{sect.reified} and \ref{exchangeability}).
\par
We begin, in Section \ref{motivating}, by setting up a motivating
example involving the Bayesian calibration of an expensive ocean model
for learning about the real ocean temperatures. In Section \ref
{meaning}, we discuss the meaning of a subjective Bayesian analysis and
develop the foundational arguments and machinery needed in order to
develop our methodology. Section \ref{sect.reified} describes our
approach in the case of finite, quantifiable, alternative modelling
judgements. Section \ref{exchangeability} extends the approach to the
case where we have infinitely many alternative modelling judgements by
using co-exchangeability to partition the space of alternatives.
Section \ref{application} illustrates the application of posterior
belief assessment to the example introduced in Section \ref
{motivating}, and Section \ref{discussion} contains discussion.

\section{A motivating example}\label{motivating}
Throughout the paper we will discuss current judgements and alternative
judgements regarding our models and specification for the prior and
likelihood in complex statistical problems, and we will claim that
perhaps none of these actually reflect our internal expert
probabilistic judgements on the observables or any key elements of the
problem for which we intend to make probabilistic inference (if such
probabilistic judgements could ever exist or be obtained if they did).
To motivate this discussion, we introduce an example from our own work
that is sufficiently complex to illustrate these ideas and to make
concrete some of our terminology. We apply our methodology to this
example in Section \ref{application}.
\par
Our example involves the calibration of a computationally expensive
ocean model used to learn about the state of the real ocean. We
describe the particular ocean model we are using and our experiments on
it in detail in Section \ref{application}. To keep things more general
initially, we will denote the ocean model $f(x, d)$ with $x$
representing the model input parameters chosen for a particular run and
$f(\cdot, d)$ representing the model output of interest, which, in this
case will be the global mean temperature at depth $d$ (in fact, the
model will output 6 hourly values for temperature, velocity and
salinity over a 3D mesh covering the globe at 31 depth levels, and we
post-process it to only work with outputs of interest).
\par
Computer model calibration \citep[developed in ][]{kennedyohagan01} is
a method for combining computer model output and real world data in
order to probabilistically describe past, current or future states of
the real world. To motivate our methodology, we describe a statistical
model for calibrating our ocean model using observed temperature at a
number of depths in order to predict temperature at an unobserved
depth. What follows is a brief outline of the statistical modelling
involved in calibration in order to set up a Bayesian network (BN) that
enables us to easily talk about required judgements and not ruled out
alternatives, and to have an example in mind when we develop the
underlying theory behind our approach. We do not reproduce the details
of the calculations required to obtain conditionals for posterior
sampling via MCMC. These are available in \cite{kennedyohagan01}.
\par
\begin{figure}[ht]
\includegraphics{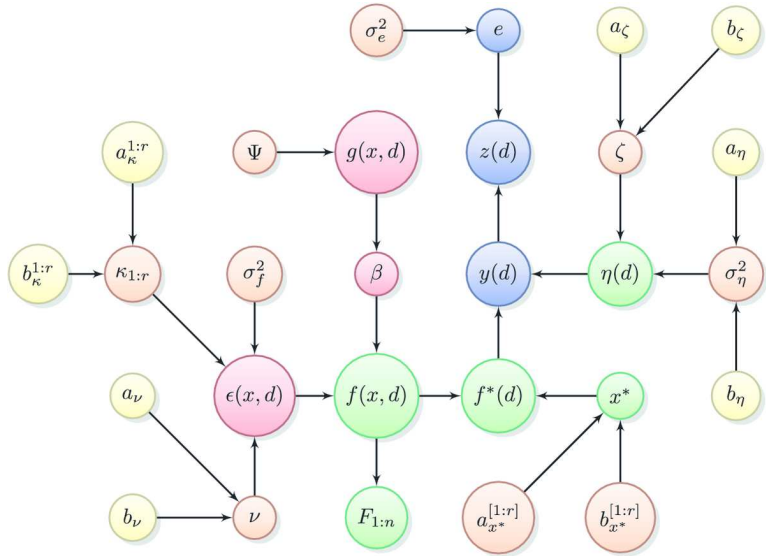}
\caption{Bayesian network for the ocean model calibration
problem.}\label{tBN}
\end{figure}

The modelling begins with observations $z(d)$ at depth $d$ of
components of the ocean $y$ made with independent error $e$ so that
\begin{displaymath}
z(d) = y(d) + e; \qquad e \sim\mathrm{N}(0 , \sigma_e^2).
\end{displaymath}
These elements relating to the real world only appear in blue on our BN
shown in Figure~\ref{tBN}.
%
\par
Components relating the ocean model to the ocean appear in green in
Figure \ref{tBN}. The model refers to the best input approach, which
indicates that there is a particular setting of the inputs, $x^*$, that
is sufficient for the ocean model in informing us about reality at any
depth $d$. The model is
\begin{displaymath}
y(d) = f(x^*,d) + \eta(d);
\end{displaymath}
with discrepancy $\eta(d)$ independent from $x^*$ and $f(x,d)$ for all
$x$. The model runs that inform us about $f(x,d)$ are denoted
$F_{1:n}$. Next we use a statistical model called an emulator to model
our uncertainty about $f(x,d)$. Components of the emulator appear as
the red nodes in Figure \ref{tBN}. The emulator is
\begin{displaymath}
f(x,d) = \beta^Tg(x,d) + \epsilon(x,d); \qquad\epsilon(x,d) \sim GP(0,
c_f(\cdot,\cdot))
\end{displaymath}
with $g(x,d)$ a chosen vector of basis functions, $\beta$ uncertain
coefficients and $\epsilon(x,d)$ an independent (from $\beta$) mean
zero Gaussian process with specified covariance function $c_f((x,d),
(x',d'))$. We usually take this to be weakly stationary so that\break
$c_f((x,d), (x',d')) = c_f(|(x,d)-(x',d')|)$.
\par
Parameters that must be specified to characterise our uncertainty for
each model component described so far appear on the BN as orange nodes.
Starting with the discrepancy, $\eta(d)$, \cite{kennedyohagan01}
suggest a mean zero weakly stationary Gaussian process prior with
Gaussian covariance function $c_{\eta}(d, d') = \sigma_{\eta}^2\exp\{
-\zeta|d-d'|^2\}$, so that a variance parameter $\sigma_{\eta}^2$ and a
depth correlation parameter $\zeta$ are introduced. The best input
$x^*$ is given a uniform distribution for each of the $r$ dimensions of
the input space with lower and upper bounds, $a_{x^*}^{[i]}$ and
$b_{x^*}^{[i]}$.
\par
We will build separate univariate emulators for each depth $d$,
altering the specification from one of requiring $\beta^Tg(x,d)$ and
$\epsilon(x,d)$ to one requiring $\beta_d^Tg_d(x)$ and $\epsilon_d(x)$.
For the response surface $\beta_d^Tg_d(x)$ we restrict the prior
specification to one of choosing the number and type of basis functions
to enter into each $g_d(x)$, choices we identify with the symbol $\Psi$
and whose options will be detailed in Section \ref{nroy.sect}. We
specify a joint prior $p(\beta, \sigma_f^2)$ as described in Section
\ref{nroy.sect}. For the Gaussian process residual, $\epsilon_d(x),$ we let
\begin{displaymath}
c(x, x') = \sigma_f^2\left(\nu\indicator_{x=x'} + (1-\nu)R(|x-x'|;\kappa
_{1:r})\right)
\end{displaymath}
where $\nu$, traditionally called the ``nugget'' \citep
{andrianakischallenor12}, represents the proportion of the residual
variability that is, in this case, due to internal variability of the
climate model \citep[see ][ for discussion]{williamsonblaker14} and
$R(|x-x'|;\kappa_{1:r})$ is a correlation function depending on
roughness parameters $\kappa$ on each dimension of the input space. We
say more on the possible choices of correlation function below.
\par
For each of the model parameters, we might be uncomfortable in
expressing values directly and could put hyper-priors on each. In our
Bayesian network we do this for $4$ of the more difficult parameters to
consider: the correlation lengths and variance of the discrepancy,
$\zeta$ and $\sigma^2_{\eta}$, the nugget term $\nu$ and the $r$
correlation parameters, $\kappa_{1:r}$, of the emulator residual. We choose
\begin{displaymath}
\zeta\sim\mathrm{G}(a_{\zeta}, b_{\zeta}), \qquad\sigma_{\eta}^2
\sim\mathrm{IG}(a_{\eta},b_{\eta}), \qquad\nu\sim\mathrm{Be}(a_{\nu
},b_{\nu}),
\end{displaymath}
and use separate $\mathrm{Be}(a_{\kappa}^{i},b_{\kappa}^{i})$ priors
for each half length correlation of the emulator residual. A half
length correlation for the $i$th input between the correlated parts of
$\epsilon(x)$ and $\epsilon(x')$ is the value of $R(|x-x'|)$ when all
elements of $x$ and $x'$ are equal with the exception of $x_i$ and
$x'_i$ and where $|x_i - x'_i|$ is equal to half of the possible range
of $x_i$ \citep[see ][ for further details]{higdonetal08,
williamsonblaker14, williamson14}. Hyper-parameters appear on our BN in yellow.
\par
The Bayesian network for this problem in Figure \ref{tBN} required many
complex probabilistic modelling judgements just to write down, and
requires a great deal more specification before a Bayesian analysis can
be performed. Values for all hyper-parameters and other parent nodes
are required to specify prior distributions, and modelling choices such
as the form of the correlation function for the Gaussian process are
required before we have established a likelihood. Even at this stage,
the required modelling and specification is, perhaps, too advanced to
expect an ocean modeller, or even a statistician-modeller team to ever
be able to complete the specification in such a way so that every
probabilistic statement we have made about any collection of nodes on
the graph represents what they really believe.
\par
In fact, as we believe would often be the case with a complex
statistical model, it is very easy to criticise the modelling we have
done as not representative of the beliefs of any expert ocean modeller.
We have stated already that in order to complete a subjective
elicitation in this problem directly from the expert, we would require
an ocean modeller with a deep understanding of probabilistic reasoning
and the implications of this type of statistical model. However, even
if such an expert did exist and had enough time to spend with us, there
are a number of criticisms they could level at the model right off the
bat. First, there is the idea that the discrepancy is additive and
independent of the model. This is problematic because if we knew $x^*$,
our ocean reconstruction $y$ would not be dynamically consistent, i.e.
it would not obey the Navier--Stokes equations over a reasonably short
time period. This is because our modelling of $\eta$ does not require
that it does (nor could it ever, reasonably). We might visit every part
of our modelling and point out holes in the realism, though we satisfy
ourselves with one further example. We know that the residual from the
response surface we fit using our emulator is not really a realisation
of a weakly stationary Gaussian process. Though, if we are lucky, we
might find a reasonable fit using such a process on a scale similar to
the average distance between our design points, it is unlikely that
this model adequately captures the behaviour of the function at, for
example, much smaller scales in all parts of the parameter space. The
assumption of stationarity here is one made largely for computational
convenience and one made because experience in this type of modelling
has often shown that breaking this assumption does not gain a great
deal in terms of predictive power relative to the effort required in
its implementation.
\par
There are many judgements already made that we might view as equivalent
to or ``just as good as'' alternative judgements that we might have
made (for example, we could have chosen a lognormal prior for $\zeta$
instead of a gamma). There are many more we are yet to make to which
the same argument will apply and, we know, for any alternative set of
judgements, that the framework we are using is only a model for our
actual underlying judgement. We make this precise in our theoretical
development below.
\par

\section{Interpreting a Bayesian analysis}\label{meaning}

There are different views as to the primitive concepts underpinning
the Bayesian approach to inference. In this paper, we are driven by
the requirement to give a formal and operational treatment which deals
with the difference between the inferences made within a traditional
Bayesian model for a complex problem and the actual inferences that we
are confident to assert on the basis of our modelling and analysis. It
is easy to see why there will often be a mismatch between these two
levels of inference. The full Bayes model requires prior specification
over complex high dimensional spaces for which we may struggle to make
meaningful judgements to the extreme level of detail required for the
analysis. Lacking confidence in our prior assessments and in the
many, often somewhat arbitrary, modelling choices that we have made in
order to complete our analysis, it is hard to have confidence in our
posterior judgments. We need to express our uncertainty about the
relevance of our Bayesian modelling to our actual real world
inferences. However, an attempt to make such a higher level
uncertainty specification through the medium of probability is, in
many ways, more difficult than was the original modelling task
itself. Partly, this difficulty is intrinsic to the task
itself. However, partly this is because probability, due to the
exhaustive level of detail that is required in the analysis, is not an
appropriate primitive, in our view, to bridge the gap between a formal Bayes
analysis and the actual posterior judgements that we may wish to
make.

Therefore, in this paper, we treat expectation, rather than
probability, as the primitive for this higher level analysis. In this
way, we may still consider any probabilities that are of direct
concern, as probabilities, in this view, are simply the expectation of
the corresponding indicator functions. However, we have the option of
restricting the level of detail of our specification, by just
considering those expectation specifications which are directly
relevant to the problem at hand. We shall show, with expectation as
primitive, how to develop a formal and operational approach to
addressing the gap between Bayesian inference and actual real world
inference for complex models, and illustrate the practicality of the
approach in application to a model of realistic size and
complexity.

It would be interesting to compare our approach with a
corresponding alternative developed purely within the formalism of
probability, but we do not know of any practical way to do this. If we
were to choose
probability as primitive, then we would need to assess a full
probabilistic specification describing the relationship between the
Bayes analysis and our judgements, the specification of which would be
even more difficult to attribute precise meaning to than the
probabilistic specification underpinning the original Bayes analysis.
\subsection{Expectation as primitive}
In the fundamental
volumes summarising his life's work \citep{definetti74, definetti75},
de Finetti based the theory of probability on expectation, or, as he termed
it, \textit{prevision}, as the primitive for the theory. de Finetti
gives the
following operational definition for your expectation, $\E{X}$, for the
random quantity $X$, namely that it is your preferred choice for the
value of $c$ when confronted with the penalty $(X-c)^2$, where the
pay-off is in some appropriate units, for example, probability
currency. (In this theory, all expectations are the subjective
judgements of specified individuals. There is no place for the notion
of objective probability, beyond, informally, the common language
meaning of a consensus judgement shared among a group of individuals.)
Of course, in practice, we may choose many alternative ways to elicit
your choice for $\E{X}$, but this definition is particularly
appropriate when deriving the formal properties of expectation as it
leads directly to the geometric structure underpinning expectation
based approaches to inference \citep[see, for example,][chapter
3]{goldsteinwooff07}.

Within this formulation of probability theory, we may derive an account
of the relationship
between conditional expectation, $\E{X|D}$ given observation of a
member of the partition $D = (D_1, \dots, D_k)$, for which one and only
one of the outcomes $D_i$ will be observed, and the actual expectation,
or prevision,
$P_{t}(X)$ that we may specify at the future time $t$ when we do observe
the outcome of the partition (We use $P_t(X)$ rather than $\Ef{t}{X}$
here, as the latter is reserved for an operation that we will define
later and make a great deal of use of. Only the currently unobtainable
time $t$ prevision is denoted with $P$. For all other subjective
expectations we use the more conventional $\mathrm{E}$ operator). This
development was first described in \cite{goldstein97} and the remainder
of this section gives a summary of it.

As the conditional probability of an event is equivalent to the
conditional expectation for the indicator function for the event, this
also gives an account of the relationship between conditional
probabilities, for example, assessed by Bayes theorem, and actual
posterior probabilities when the conditioning event is observed. It is
important to maintain the distinction between conditional
probabilities, defined in subjective probability theory as bets which
are called off unless the
conditioning event occurs, and posterior probabilities which you
assign when you have seen the corresponding event. There is no obvious
formal relationship between these two notions whatsoever, and
foundational descriptions typically restrict the formal development to
perfectly rational individuals, operating in small and very tightly
constrained worlds. As such, conditional probability may be
interpreted as a simple and flexible model for real world
inference. However, like any model, it is essential to be careful in
considering the relationship between the model and the real world,
which is not tightly constrained and concerns the actual inferences of
real individuals.

Although there is no deterministic relationship between conditional
and posterior probabilities, we can derive certain probabilistic
relationships between the two concepts. Specifically, we may use the
notion of temporal sure preference. Suppose that we must choose between
two random
penalties $A$ and $B$. We say that we have
a sure preference for $A$ over $B$ at (future) time $t$ if we are sure
now that, at time $t$, we will prefer $A$ to $B$. The temporal sure
preference (TSP) principle says that if we have a sure preference for
$A$ over $B$ at $t$, then we should not have a preference for $B$ over
$A$ now. Temporal sure preference is a very weak property (as it is
hard for us to hold temporal sure preferences). However, it is
sufficient to derive the basic probabilistic relationship between
conditional and posterior expectations, which is as follows.

Suppose that we currently specify our conditional expectation for a
random vector $X$ given partition $D=(D_1, \dots, D_k)$. Suppose, at
future time $t$ when we have observed which element of $D$ occurs, we
make an actual posterior expectation statement $P_{t}(X)$. If we denote
by $\E{X|D}$ the random quantity which takes value $\E{X|D_i}$ if we
observe outcome $D_i$, then given TSP we can show \citep[see][for
details]{goldstein97} that, now, we must
have the following orthogonal decomposition:
%
\begin{equation}
\label{TSPrel}
X = (X - P_{t}(X)) \oplus(P_{t}(X) - \E{X|D}) \oplus(\E{X|D} - \E{X})
\oplus\E{X}.
\end{equation}

Each of the three bracketed terms on the right hand side of
(\ref{TSPrel}) has expectation zero, and the three terms are mutually
uncorrelated (the notation $A \oplus B$ means that the random vectors
$A$ and $B$ are uncorrelated). More strongly, $(P_{t}(X) - \E{X|D})$ has
expectation zero conditional on each member of the partition. We may
view $\E{X|D}$ as a prior inference for $P_{t}(X)$, in the sense that
$\E{X|D}$ bears the same relationship to $P_{t}(X)$ that $P_{t}(X)$ bears
to $X$. In this way, we can make precise the notion that the formal
Bayes analysis is a model for your posterior judgements, whose
relationship is the same as for any model of a real world process \citep
{goldstein11}.

As the terms in (\ref{TSPrel}) are uncorrelated, we can interpret this
construction in terms of the resulting variance partition,
%
\begin{equation}
\label{VarTSPrel}
\var{X} = \var{X - P_{t}(X)} + \var{P_{t}(X) - \E{X|D}} + \var{\E{X|D}}.
\end{equation}

The variance of $X$ is decomposed into three components. $\var{\E{X|D}}$
expresses the information content of the formal Bayes analysis,
$\var{P_{t}(X) - \E{X|D}}$ expresses the information content of the further
judgements and information, in addition to the formal Bayes analysis,
that we may bring to bear by time $t$, and $\var{X - P_{t}(X)}$ expresses
the intrinsic limitations to our inferences by this time. We can
expand the value of the Bayesian component of the analysis by
expanding the partition and thus the variation attributable to the
conditioning.

While we have described the variance partition for a full Bayes
analysis, the argument from TSP actually establishes the more general
property that for any random quantity $W$ which we will certainly
observe by time $t$, then we must assign
\begin{displaymath}
\E{(X_i - P_{t}(X_i))^2} \leq\E{(X_i-W)^2}.
\end{displaymath}
By varying the choice of $W$, we can derive the corresponding representation
%
\begin{equation}
\label{TSPproj}
X = (X - P_{t}(X)) \oplus(P_{t}(X) - \Ef{U}{X}) \oplus(\Ef{U}{X} - \E
{X}) \oplus\E{X}
\end{equation}
where $\Ef{U}{X}$ is the orthogonal projection of $X$ into the collection
of linear combinations of the elements of the collection $U = (W_1,
\dots, W_k)$ where each $W_i$ is a random quantity which will be
observed by time $t$, and the terms on the right hand side of
(\ref{TSPproj}) obey the same orthogonality properties as those of
(\ref{TSPrel}), with a similar interpretation.

In this next section, we shall combine these two variance partitions to give
a general representation of the inference that we may make through a
formal Bayes analysis for a complex problem. Firstly, (\ref{TSPrel})
is based on the requirement that all of the conditioning statements
refer to your actual probabilistic judgements. In practice, for
complex problems, the level of detail required is too extreme for this
to be possible. We should instead view the choices made in the full
Bayes analysis as providing a model for your actual probabilistic
judgements (in a way we shall make clear in the next section), which is
itself a model for your posterior
inferences. Therefore, we will consider a range of modelling choices
for likelihoods, priors and so forth. Each choice leads to a Bayesian
conditional expectation. We now use (\ref{TSPproj}) to construct the
corresponding relationship between the outcomes of each collection of
modelling choices that have been considered and the actual inferences
that you are able to make given the data and your collection of
analyses.

\section{Posterior belief assessment}\label{sect.reified}
Define a vector of quantities of interest, $\bel$, so that your
principal interest in your Bayesian analysis is in finding posterior
expectations of $\bel$ given relevant data $\dat$. Note that this
definition allows us to be interested in posterior variances and
posterior probabilities of key events (defined as expectations of
indicator functions), for example, the probability that some quantity
exceeds a threshold. Define your current modelling judgements as the
set $\jdg{0}$. These judgements comprise everything required to perform
the Bayesian calculation, including form and any hyper-parameters of
both the prior and the likelihood. We do not believe all of the
probability statements in $\jdg{0}$, but we do view $\jdg{0}$ as
``representative enough'' of the structure of the problem and any
scientific judgements and beliefs that we do hold, to allow us to
perform the Bayesian computation and to believe that $\E{\bel|\dat;\jdg
{0}}$ is informative for our posterior prevision for~$\bel$. As
discussed in Section \ref{meaning}, this conditional expectation will
not be your prevision. Operationally, this means that we prefer penalty
$(\bel_i - P_t(\bel_i))^2$ to $(\bel_i - \E{\bel_i|\dat;\jdg{0}})^2$,
equivalently
\begin{displaymath}
\E{(\bel_i - P_t(\bel_i))^2} \leq\E{(\bel_i - \E{\bel_i|\dat;\jdg{0}})^2}.
\end{displaymath}
\par
In Section \ref{meaning}, we argued that nothing we can do prior to
time $t$ will give us $P_t(\bel)$ at time $t$ after seeing $z$. In
particular, we cannot set up a called off bet via our Bayesian analysis
under $\jdg{0}$ and be required to hold $\E{y|z;\jdg{0}}$ at time $t$,
even if we really believed all of the probability statements in $\jdg
{0}$. However, though $P_t(\bel)$ may be unobtainable within the
constraints of the Bayesian formalism that we have expressed, we would
like to get as close as possible to it using the fully Bayesian
machinery. For example, can we find good choices of $\Gamma(\dat)$ to
improve upon $\E{\bel|\dat;\jdg{0}}$ (in the sense that you prefer
penalty $(\bel_i - \Gamma_i(\dat))^2$ to $(\bel_i - \E{\bel_i|\dat;\jdg
{0}})^2$ for any vector-valued function $\Gamma$ of the data at time
$t$, for all $i$)?
\par
Though the $\jdg{0}$ contains sufficient structure to enable us to
perform a Bayesian analysis that we view as informative for certain
posterior expectations over quantities of interest, as we do not
believe all of the probability statements made by $\jdg{0}$, we do not
view these judgements as uniquely representative of the problem
structure and our scientific insight. By changing parts of the model
and/or the prior (or even the way we generate posterior samples) in a
way that either represents our scientific insights and beliefs in a
different way, or that changes them in such a subtle way that we would
still view the analysis as informative for our prevision for $\bel$, we
can arrive at a set of alternative judgements $\jdg{1}, \ldots, \jdg
{k}$. We could, in principle, run alternative Bayesian analyses under
each set of judgements and compute $\E{\bel|\dat;\jdg{i}}$ for $i= 1,
\ldots, k$. We assume a finite collection of alternative judgements we
would adopt here, and relax this constraint in Section \ref{exchangeability}.
\par
We note here that in viewing the Bayesian analysis like this, we are
using probability in two different ways. The first, in the true
subjective sense, represents what we actually believe about a quantity,
and these are our previsions. The second uses probability as a useful
modelling language for transforming our judgements into previsions. The
probability statements made as part of any $\jdg{i}$ are therefore not
belief statements, nor do we assume that there is an underlying true
representation of our beliefs, $\jdg{*}$ say, that could be elicited if
we could think hard enough for long enough to extract it. However, it
is reasonable to talk of potential sets of judgements not representing
our beliefs well enough, so that we could rule them out as being
informative for our previsions. For this reason, we will sometimes
refer to the set of alternative judgements $\jdg{1}, \ldots, \jdg{k}$
as not ruled out judgements.
\par
That we have chosen $\jdg{0}$ over any $\jdg{1}, \ldots, \jdg{k}$, may
lead us to the view that we prefer penalty $(\bel_i - \E{\bel_i|\dat
;\jdg{0}})^2$ to $(\bel_i - \E{\bel_i|\dat;\jdg{i}})^2$ for $i=1, \ldots
, k.$ However, some of the elements of $\jdg{0}$ may have been chosen
for pragmatic reasons. For example, we might choose a conjugate prior
family to reduce the computational burden in part of our posterior
sampling scheme. Even if we were to adopt the view that $\jdg{0}$ leads
to conditional expectations that we do prefer to conditional
expectations derived through Bayesian analyses under any of the
alternative sets of judgements, there may be derivable random penalties
that make use of the collection of our current and alternative
judgements that we prefer over any based on one set of them. We now
consider combinations of quantities related to our alternative judgements,
$\jdg{1},\ldots,\jdg{k}$. We show that there is a linear
combination of the collection of conditional expectations that is at
least as close to your actual prevision as $\E{\bel|\dat;\jdg{0}}$ and
derive its properties. Define $\G$ to be the vector
($\E{\bel|\dat;\jdg{0}}$, $\E{\bel|\dat;\jdg{1}}$,$ \ldots$,
$\E{\bel|\dat;\jdg{k}}) = (\G_1, \ldots, \G_{k+1})$. Let $\G_0$ be the
unit constant. We have the following result.
\begin{thm}
Let
%
\begin{equation}\label{reified}
\Ef{\G}{\bel} = \E{\bel} + \cov{\bel}{\G}\var{\G}^{-1}(\G- \E{\G}),
\end{equation}
an expectation which we term, our posterior belief assessment for $\bel
$. Then
\begin{itemize}
\item[(i)] $\Ef{\G}{\bel}$ is at least as close to $\bel$ as $\E{\bel
|\dat;\jdg{0}}$. Equivalently, for each $i$,
\begin{displaymath}
\E{(\bel_i - \Ef{\G}{\bel_i})^2} \leq\E{(\bel_i - \E{\bel_i|\dat;\jdg{0}})^2}.
\end{displaymath}
where $\Ef{\G}{\bel_i}$ is the $i$th component of $\Ef{\G}{\bel}$.
\item[(ii)] $\Ef{\G}{\bel}$ is at least as close to $P_t(\bel)$ as $\E
{\bel|\dat;\jdg{0}}$. Equivalently, for each $i$,
\begin{displaymath}
\E{(P_t(\bel_i) - \Ef{\G}{\bel_i})^2} \leq\E{(P_t(\bel_i) - \E{\bel
_i|\dat;\jdg{0}})^2}.
\end{displaymath}
\end{itemize}
\end{thm}
\begin{proof}
Our posterior belief assessment is the Bayes linear rule for $\bel$, given
$\G$ \citep[see][]{goldsteinwooff07}. For each $i$, $\Ef{\G}{\bel_i}$ is
the Bayes linear rule for $\bel_i$, namely the linear combination $\sum
_{i=0}\alpha_i\G_i$ that
minimises
\begin{displaymath}
\E{\left(\bel_i -
\sum_{i=0}\alpha_i\G_i\right)^2}.
\end{displaymath}
\citep[see page 56 in][for a proof of this
result]{goldsteinwooff07}. Therefore, property (i) follows immediately.
\par
From TSP, $\bel- P_t(\bel)$ has mean zero and is uncorrelated
with all random quantities that will be known by time $t$, and, in
particular, all the elements of $\G$. Therefore, for each $i$,
\begin{align*}
\E{(\bel_i - \Ef{\G}{\bel_i})^2} &\leq\E{(\bel_i -
\E{\bel_i|\dat;\jdg{0}})^2} \quad \Rightarrow\\
\E{(\bel_i -P_t(\bel_i))^2} + \E{(P_t(\bel_i) - \Ef{\G}{\bel_i})^2}
&\leq \E{(\bel _i - P_t(\bel_i))^2}\\
&\quad + \E{ (P_t(\bel_i) - \E{\bel_i|\dat;\jdg{0}})^2} ,
\end{align*}
and property (2) follows immediately.
\end{proof}
Our theorem identifies the key quantities that we need to assess in
order to implement the method, namely $\E{\G}$, $\var{\G}$ and $\cov
{\bel}{\G}$. In Section \ref{algorithm}, we present a method for
evaluating these quantities.
Note that the theorem implies that by computing the posterior belief
assessment using additional experiments, we resolve an additional
proportion of our uncertainty in $\bel$, compared with the case where
we only have one Bayesian analysis. We can estimate a lower bound on
the proportion of uncertainty reduced by observing the ratio of the
adjusted variances $\var{\bel_i - \Ef{\G_1}{\bel_i}}$ and $\var{\bel_i
- \Ef{\G}{\bel_i}}$ (note $\Ef{\G_1}{\bel_i}$ is computed using (4)
replacing $\G$ with $\G_1$). This is a direct corollary of part (i).
\par
Readers familiar with Bayesian Model Averaging \citep
[BMA,][]{hoetingetal99}, may notice that (\ref{reified}), which is
simply a constant plus a linear combination of each alternative
posterior expectation, may be interpreted as a constant plus the BMA
posterior expectation under a posterior distribution over model $\jdg
{i}$ equal to $W_{i+1}$ with
\begin{displaymath}
W_j = (\cov{\bel}{\G}\var{\G}^{-1})_j.
\end{displaymath}
However, our approach is not the same as we have made expectation
primitive and not probability. In BMA, one needs a probability
distribution over all possible models in order to proceed. Hence, in
this context, the BMA cure to needing to make so many judgements in a
full Bayesian analysis and to not being sure that any of these
judgements represent our actual beliefs, is to make infinitely more
judgements (in the form of a probability distribution across all
possible models and conditional probabilities under each model), and
proceed from there. But now we are back in the same boat in which we
started, requiring many more prior judgements in order to do the
Bayesian analysis. We can ask the same questions regarding whether we
actually hold each of these judgements, whether the statistical
modelling is representative of them, or if there are alternatives that
we feel should be explored and so on. More subtly, as mentioned in the
introduction, BMA must also assume that one of our alternative models
is the truth, and we are unwilling to make this assumption.\vadjust{\eject}
\par
With expectation as primitive in a posterior belief assessment, we
avoid the infinite regress. We make no probabilistic prior statements
about any of our alternative models, hence do not assume any represent
our true underlying beliefs, and simply seek to get closer to actual
prevision, $P_t(\bel)$, than $\E{\bel|\dat;\jdg{0}}$ via orthogonal
projection of $\bel$ onto a vector $\G$ of alternative Bayesian
analyses. Importantly we have established the relationship between our
posterior belief assessment, $\Ef{\G}{\bel}$, and what we will actually
think at time $t$, $P_t(\bel)$, which would be a task still ahead of an
analyst adopting BMA.
\subsection{Practical considerations}\label{practical}
The above account states that if we can write down a finite set of
alternative judgements $\jdg{1}, \ldots, \jdg{k}$, to $\jdg{0}$, and if
we can obtain conditional expectations for key quantities, $\bel$,
given data $\dat$, using a full Bayesian analysis under each $\jdg{i}$,
$i=0, \ldots, k$, then we can perform a posterior belief assessment to
get closer to our actual judgements if we can compute (\ref{reified}).
However, there are a number of practical considerations to doing this,
even if we can write down $\jdg{1}, \ldots, \jdg{k}$.
\par
One such consideration is that each alternative Bayesian analysis must
be performed in order to obtain each $\E{\bel|\dat;\jdg{i}}, i =
1,\ldots,k,$ so that we may evaluate $\G$. There are, essentially, two
factors impacting upon the feasibility of this step. The first involves
the availability of extra computer power. For many of the $\jdg{i}$, we
may only have small changes to hyper-prior quantities to make, or
alternative, yet trivially computable, formulations of the likelihood.
In these cases, performing the alternative Bayesian computations (e.g.
via MCMC) involves little or no manpower and instead requires more
computing. In these cases the alternative analyses can run in parallel
either on single multi-core machines, or via clusters or distributed
computing, representing little additional complexity to the overall
analysis. The second factor to consider regards the additional
complexity of the Bayesian computation for any of the alternative
judgements $\jdg{i}$. An alternative model may lead to a Bayesian
calculation that we don't yet know how to do well, or that would take
so much additional effort to implement that we regard performing this
alternative analysis in order to provide a posterior belief assessment
to be either impractical (we do not judge it to be a good use of our
time), or infeasible given our current budget constraints. We visit a
solution to this issue in Section~\ref{exchangeability}.
\par
An important consideration is that we require a number of quantities in
order to compute $\Ef{\G}{\bel}$ via (\ref{reified}). In particular, we
must either specify or compute $\E{\bel}$, $\E{\G}$, $\var{\G}$ and
$\cov{\bel}{\G}$. Prior to offering practical suggestions for how this
might be achieved, we note that our account thus far has presented an
argument for why (\ref{reified}), if it can be computed, is closer to
actual prevision and is thus important. If the meaning of any posterior
statements that we will make about $\bel$ having obtained data $\dat$
is important, we have provided a route towards clear meaning through
(\ref{reified}). Specifically, we have shown that $\Ef{\G}{\bel}$ is as
close as possible with respect to squared error loss to time $t$
prevision of the true $\bel$. This represents a novel foundational and
methodological step, focusing our attention on the further key
quantities we require to claim ownership of posterior belief
statements. As such, an interesting avenue for further research in this
matter could focus on precisely how key quantities for posterior belief
assessment should be obtained.
\par
That said, we now discuss two particular ideas. The first is to focus
elicitation on $\E{\bel}$, $\E{\G}$, $\var{\G}$ and $\cov{\bel}{\G}$.
The elicitation of purely second order judgements is necessarily
simpler than the sort of careful elicitation required for a fully
subjective Bayesian analysis in a complex problem. This is discussed in
many sources, an overview is in Chapter 2 of \cite{goldsteinwooff07}.
For example, though an expert may find it extremely hard to think about
distributions for hyper-parameters of distributions on other parameters
that purportedly control the distribution of some quantity related to
the complex system $\bel$ (wherein her expertise lies) through a
conditional probability model, she will be far more comfortable
considering $\E{\bel}$ and $\var{\bel}$ as it involves thinking about
the quantity she understands directly. Though simpler than a fully
probabilistic elicitation exercise for the statistical model in
question, the above argument may not be viewed, in certain problems, to
be as strong when it comes to determining those key judgements
involving $\G$. Hence we provide a second suggestion involving sampling.
\subsection{An algorithm for computing posterior belief
assessments}\label{algorithm}
Our algorithm for obtaining $\E{\G}$, $\var{\G}$ and $\cov{\bel}{\G}$,
the key quantities required to compute (\ref{reified}), requires that
$\E{\bel}$ and $\var{\bel}$ are specified, either through elicitation
or otherwise. We can then form a distribution $(\E{\bel}, \var{\bel})$
and sample values of $\bel$. We also use the same arguments to form a
distribution for the observables using the mean and variance $(\E{\dat
|\bel}, \var{\dat|\bel})$. For example, $\dat$ might be a measurement
on $\bel$ made with mean zero, uncorrelated error with given variance
$\sigma_e^2$ and $\E{\dat|\bel} = \bel$ and $\var{\dat|\bel} = \sigma
_e^2$. We can now sample a value of $\bel$ from $(\E{\bel}, \var{\bel
})$ and corresponding values of $\dat$ using $(\E{\dat|\bel}, \var{\dat
|\bel})$. Once we are able to generate such samples, our algorithm
proceeds as follows.
\setlength\leftmargini{38pt}
\begin{itemize}
\item[\textbf{Step 1}] Sample a value $\hat{\bel}$ and $\hat{\dat}$
from $(\E{\bel}, \var{\bel})$ and
$(\E{\dat|\hat{\bel}}, \var{\dat|\hat{\bel}}),$ respectively.
\item[\textbf{Step 2}] Use the full Bayes machinery to compute each $\E
{\bel|\hat{\dat};\jdg{0}}$, $\E{\bel|\hat{\dat};\jdg{1}}, \ldots$, and
use them to form $\hat{\G}$.
\item[\textbf{Step 3}] Repeat this process to obtain a large
    number, $N$, of sample pairs $(\hat{\bel}_1, \hat{\G}_1),
    \ldots,\break (\hat{\bel }_N, \hat{\G}_N)$.
\item[\textbf{Step 4}] Assess $\E{\G}$, $\var{\G}$ and $\cov{\bel}{\G}$
by computing the sample means and variances of the $\hat{\G}s$ and
their covariance with the $\hat{\bel}$s.
\end{itemize}
\par
To make the sampling algorithm computationally feasible using
distributed computing, we may have to, for example, reduce the number
of samples in our MCMC at step~2. When we do not have a firm judgement
as to our choice for the distribution of $\bel$ or $\dat|\bel$, then,
for negligible computational cost, we may repeat the analysis under a
variety of choices, by re-weighting the samples that we obtain using
our algorithm, under each such choice.
\par
Performing such alternative analyses by re-weighting may be important
in the case where we have no firm choice for these distributions or
where the model we use to obtain $\dat|\bel$ is tentative or one of a
collection of models we are comfortable with. Depending on the problem,
there is a potential for the impact on the posterior belief assessment
of these choices to be non-negligible. If, through re-weighting to
other reasonable forms of distribution for $\bel$ and $\dat|\bel$, we
establish this to be the case, we would invest more time and effort
eliciting features (such as skew, kurtosis, certain percentiles of the
distribution function) of these distributions from our experts. As
these quantities (each of which is an expectation) are functions of the
observables, we view this as easier for a domain expert than
elicitation of parameter distributions in our complex statistical
models. However, we do not require that any elicited distribution here
represent the subjective beliefs of the analyst or expert in the sense
that every or any probability statement made is believed, as we are
using probability here in the way described near the beginning of
Section \ref{sect.reified}: as a tool to assist with obtaining a
handful of beliefs we are prepared to adopt as prevision through the
machinery of posterior belief assessment. The impact of these choices
on the results of our sampling algorithm and the posterior belief
assessment itself, as well as formal methods for assessing this and
handling large impacts would be an interesting avenue for future
research in this area.

\subsection{Comparison with robust Bayes and sensitivity
analysis}\label{sect.robust}
A robust Bayesian analysis requires us to be able to explore,
analytically (usually), the behaviour of any key posterior quantities
under all possible classes of prior and likelihood combination of
interest. This is an extremely complex and arduous task, and may be
infeasible for complex models. Even if completed, however, unless a
particular posterior quantity is (practically) invariant to all
choices, it is unclear what belief statements you can make and how they
might relate to your prevision. Similarly, though a sensitivity
analysis has the potential to be just as thorough as a posterior belief
assessment in terms of the alternative models and priors explored ($\jdg
{1}, \ldots, \jdg{k}$ in our account thus far), it leaves any formal
relationship between the quantities computed under alternative
judgements, your beliefs and prevision unestablished. Typically, we
might report posterior beliefs under $\jdg{0}$ and use the alternative
analyses to comment on sensitivity of these beliefs to our judgements,
without actually incorporating data from any further analyses into our
beliefs given the data.
\par
A posterior belief assessment offers a practical and tractable
investigation into the robustness of any conclusions you may wish to
draw to alternative prior and likelihood choices. Though requiring
extra computing power, the calculations required are essentially
repeats of previous analyses and are therefore relatively
straightforward to undertake, regardless of the complexity of the
statistical modelling. Further, our approach reduces your uncertainty
about the quantities of interest and gives an actual quantification of
the current state of your beliefs about each quantity that accounts for
all of the additional information contained in further analyses,
through $\Ef{\G}{\bel}$.
\par
The account in this section has described posterior belief assessment
for the situation where we can write down and construct full Bayes
analyses for a finite set of alternative judgements $\jdg{0}, \jdg{1},
\ldots, \jdg{k}$. However, in many practical cases of interest, it is
unlikely that $k$ will be small enough to permit an exhaustive
analysis, or perhaps the number of not ruled out judgements is, in
principle, infinite (for example, we may have had to specify continuous
hyper parameters). In the next section, we outline how exchangeability
may be imposed on our judgements and exploited to provide posterior
belief assessment in more realistic and general situations.

\section{Co-exchangeable classes and posterior belief assessment}\label
{exchangeability}
One of the key tools in the subjective Bayesian's armoury when faced
with a collection of similar random quantities is the judgement of
exchangeability. Suppose we have collection of separate flips of a
coin, $X_1, X_2, \ldots,$ and we want to predict the outcome of some
future flip of the same coin. The judgement of exchangeability,
effectively that our judgements for any combination of the $X$s does
not depend on the indices, allows us to treat each $X_i$ as a draw from
some probability distribution, to learn about the parameters of that
distribution using our collection of flips, and to make inference about
future, yet unobserved $X$s. This is a consequence of the
representation theorem \citep{definetti74}, and a detailed introductory
discussion is given in \cite{goldstein12}.
\par
Having a collection of alternative judgements $\jdg{0}, \jdg{1}, \ldots
$ leading, in principle, to quantities $\E{\bel|\dat;\jdg{0}}$, $\E{\bel
|\dat;\jdg{1}}, \ldots$ leads us to ask questions such as are there any
exchangeability judgements we might be willing to make and what might
they be? There is no reason to assume, in general, that the collection
$\E{\bel|\dat;\jdg{0}}, \E{\bel|\dat;\jdg{1}}, \ldots$ should be
exchangeable. For example, the way $\jdg{i}$ affects either the
modelling or the prior choices might be quite different from another
$\jdg{j}$, and this might lead us to strong a priori views on the
differences between $\E{\bel|\dat;\jdg{i}}$ and $\E{\bel|\dat;\jdg{j}}$.
\par
The key idea in this section will be to group our alternative
judgements into classes within which we are prepared to impose some
form of exchangeability on $\E{\bel|\dat;\jdg{i_1}},\break
\E{\bel|\dat;\jdg {i_2}}, \ldots$ for $\jdg{i_1}, \jdg{i_2}, \ldots$ in
class $i$. This will allow us to perform a handful of alternative
Bayesian analyses for some of our not ruled out judgements and to use
the information from these, to learn about other analyses, for example,
those that we do not have time or the ability to perform. Before
detailing our approach, we must first set up the type of
exchangeability we require, namely \textit {second order
exchangeability}.
\subsection{Second order exchangeability}
In order to establish an inferential framework within which alternative
Bayesian analyses may be used to discover relationships between
explored judgements and the rest of our not ruled out judgements, we
require further statistical modelling. We prefer to minimise the number
of additional judgements required in order to do this, and so turn to
second order exchangeability (SOE). A collection of quantities $X_1,
X_2, X_3, \ldots$ are SOE if
\begin{displaymath}
\E{X_i} = \mu; \qquad\var{X_i} = \Sigma\quad\forall i; \qquad\cov
{X_i}{X_j} = \Gamma\quad\forall i \neq j.
\end{displaymath}
In words, SOE represents an invariance of our judgements about means,
variances and pairwise covariances in a collection to labelling. For an
infinite SOE collection, the representation theorem \citep[][chapter
6]{goldsteinwooff07} gives an orthogonal decomposition of each $X_j$
into a common term $\M{X}$ and a residual $\R_j(X)$:
%
\begin{equation}\label{soe1}
X_j = \M{X} + \R_j(X)
\end{equation}
with
%
\begin{equation}\label{soeM}
\E{\M{X}} = \mu; \qquad\var{\M{X}} = \Lambda;
\end{equation}
%
\begin{equation}\label{soeR}
\E{\R_j(X)} = 0; \qquad\var{\R_j(X)} = \Sigma- \Lambda;
\end{equation}
and
%
\begin{equation}\label{soeCovs}
\cov{\M{X}}{\R_j(X)} = \cov{\R_j(X)}{\R_k(X)} = 0 \qquad\forall\,j,
\text{and}\, k \neq j.
\end{equation}
\par
The inferential power of this judgement is that it leads to the
orthogonality of a set of $n$ observations $D_n = (X_1, X_2, \ldots,
X_n)$ and any $X_j$, $j>n$, given $\M{X}$. This orthogonality, termed
\textit{Bayes linear sufficiency}, means that the Bayes linear update
$\Ef{D_n\cup\M{X}}{X_j}$, is equivalent to $\Ef{\M{X}}{X_j}$, which
implies \citep[\xch{cf.}{c.f.}][page 195]{goldsteinwooff07}
%
\begin{equation}\label{sufficiency}
\Ef{D_n}{X_j} = \Ef{D_n}{\M{X}}.
\end{equation}
Hence, suppose we judged that $\E{\bel|\dat;\jdg{1}}$, $\E{\bel|\dat
;\jdg{2}}, \ldots$, were an infinite SOE sequence and that we only
observed $D_k = \{\E{\bel|\dat;\jdg{k}}; k=1,\ldots,4\}$, then any $\Ef
{D_k}{\E{\bel|\dat;\jdg{i}}}$, $i>4$, is equal to $\Ef{D_k}{\Mno}$ and
can be computed using (\ref{soe1}), (\ref{soeM}), (\ref{soeR}) and (\ref
{soeCovs}). (Note that a Bayes linear update takes the same form as
\xch{(\ref{reified}).)}{(\ref{reified})).} \par In general, we will not wish to impose SOE
across the collection of posterior expectations under all alternative
judgements we might make, as we may have relatively strong judgements
regarding the relationship between certain alternative modelling
choices. For example, if $\jdg{i}$ and $\jdg{j}$ differ only in the
values of the hyper-parameters of one particular prior distribution,
and $\jdg{k}$ represents an alternative formulation of the likelihood,
then we might view, a priori, that
\begin{displaymath}
\cov{\E{\bel|\dat;\jdg{i}}}{\E{\bel|\dat;\jdg{j}}} \neq\cov{\E{\bel
|\dat;\jdg{i}}}{\E{\bel|\dat;\jdg{k}}}
\end{displaymath}
and similarly for $\cov{\E{\bel|\dat;\jdg{j}}}{\E{\bel|\dat;\jdg{k}}}$.
The potential existence of such relationships between posterior
conditional expectations under different judgements motivates a
partitioning of our judgements into co-exchangeable classes.
\subsection{Co-exchangeable classes}
Let $Y_1, Y_2, \ldots$ be a collection of sequences. We say that they
are \textit{co-exchangeable} if, for any fixed $k$, the sequence
$Y_{k1}, Y_{k2}, \ldots$ is infinite second order exchangeable and if,
for $j \neq k$, $\cov{Y_{jm}}{Y_{kn}} = \Sigma_{jk}$ for any $m$, $n$.
Let our alternative judgements be partitioned into $k$ classes $\jdg
{1}, \ldots, \jdg{k}$ with $\jdg{i} = \{\jdg{i_1}, \jdg{i_2}, \ldots\}
$, each chosen so that we may define corresponding co-exchangeable
classes $\class{1}, \ldots, \class{k}$ with $\class{i} = \{\E{\bel|\dat
;\jdg{i_1}}, \E{\bel|\dat;\jdg{i_2}}, \ldots\}$. Then, the
representation theorem gives
%
\begin{equation}\label{coex1}
\E{\bel|\dat;\jdg{i_j}} = \M{\class{i}} + \R_j(\class{i}).
\end{equation}
Equation (\ref{coex1}) and co-exchangeability between classes means that
%
\begin{equation}\label{coexcov}
\cov{\E{\bel|\dat;\jdg{i_j}}}{\E{\bel|\dat;\jdg{m_l}}} = \cov{\M{\class
{i}}}{\M{\class{m}}} \qquad i \neq m,
\end{equation}
and, from (\ref{soeM}), (\ref{soeR}) and (\ref{soeCovs}), we have
%
\begin{equation}\label{coexcov2}
\cov{\E{\bel|\dat;\jdg{i_j}}}{\E{\bel|\dat;\jdg{i_l}}} = \var{\M{\class
{i}}} \qquad j \neq l.
\end{equation}
\par
By partitioning into co-exchangeable classes in this way, we are able
to run a small number of alternative Bayesian analyses corresponding to
a subset of alternative judgements, to use these to update beliefs
about any other Bayesian analysis we might have performed and to
compute our posterior belief assessment. Suppose we were to compute
posterior expectations based on Bayesian analyses on $N$ sets of
alternative judgments. Let $n(1),\ldots, n(N) \in\{1, \ldots, k\}$ and
$m(1),\ldots,m(n) \in\mathbb{Z}$ be indices corresponding to the chosen
judgements, so that the $j$th Bayesian analysis performed happens under
$\jdg{n(j)_{m(j)}}$. Define \xch{$D = \{\mathrm{E}[\bel|\dat;\jdg
{n(1)_{m(1)}}],\ldots,\mathrm{E}[\bel|\dat;\jdg{n(N)_{m(N)}}]\}$}{$D =
\{\E{\bel|\dat;\jdg
{n(1)_{m(1)}}},\ldots,\E{\bel|\dat;\jdg{n(N)_{m(N)}}}\}$}, then (\ref
{coexcov}) and (\ref{coexcov2}) lead to a version of (\ref
{sufficiency}), namely
%
\begin{equation}\label{sufficiency2}
\Ef{D}{\E{\bel|\dat;\jdg{i_j}}} = \Ef{D}{\M{\class{i}}}.
\end{equation}
Further, each $\Ef{D}{\M{\class{i}}}$, for $i=1,\ldots,k$, is linear
combination of the elements of~$D$. Hence, when looking to use the
elements of $D$ to perform a posterior belief assessment, we can
simultaneously include $D$ and all of the information contained therein
regarding unobserved alternative Bayesian analyses by selecting $\G$ to
be the vector $(\E{\bel|\dat;\jdg{0}}, \Ef{D}{\M{\class{1}}}, \ldots,
\Ef{D}{\M{\class{k}}})$ $= (\G_1, \ldots, \G_k)$. Letting $\G_0$ be the
unit constant, then the results in Theorem 1 hold.
\subsection{Special case: co-exchangeability of $y$}
A further co-exchangeability judgement simplifies the task of prior
specification required in order to compute $\Ef{\G}{\bel}$: namely,
that of co-exchangeability between $\bel$ and the members of each
class, $\class{1}, \ldots, \class{k}$.
\par
Determining $\cov{\bel}{\G}$ requires determining $\cov{\bel}{\E{\bel
|\dat;\jdg{0}}}$ and\break $\cov{\bel}{\Ef{D}{\M{\class{i}}}}$ for
$i=1, \ldots, k$. Supposing we have the quantities required to compute
each $\Ef{D}{\M{\class{i}}}$, then the co-exchangeability of $\bel$ and
the members of each class $\class{1}, \ldots, \class{k}$ means that we
require only $\cov{\bel}{\E{\bel|\dat;\jdg{0}}}$ and $\cov{\bel}{\M
{\class{i}}}$ for $i=1, \ldots, k$ in order to compute
$\cov{\bel}{\G}$.
\par
To show this, we first note that
\begin{displaymath}
\begin{split}
\cov{\bel}{\Ef{D}{\M{\class{i}}}} &= \cov{\bel}{\cov{\M{\class
{i}}}{D}\var{D}^{-1}D} \\
&= \cov{\bel}{D}W_i
\end{split}
\end{displaymath}
with $W_i = \var{D}^{-1}\cov{D}{\M{\class{i}}}$. $W_i$ is required in
order to compute $\Ef{D}{\M{\class{i}}}$, hence, the additional burden
in prior specification in order to compute\break
$\cov{\bel}{\Ef{D}{\M{\class {i}}}}$ requires us to compute
$\cov{\bel}{D}$. From (\ref{coex1}),
\begin{displaymath}
D_j = \M{\class{n(j)}} + \R_{m(j)}(\class{n(j)}).
\end{displaymath}
Defining $\tilde{\Mno} = \M{\class{n(1)}}, \ldots, \M{\class{n(N)}}$,
then co-exchangeability of $\bel$ implies that
\begin{displaymath}
\cov{\bel}{D}W_i = \cov{\bel}{\tilde{\Mno}}W_i
\end{displaymath}
so that we need only specify $\cov{\bel}{\M{\class{i}}}$ for $i=1,
\ldots, k$ in order to compute $\cov{\bel}{\Ef{D}{\M{\class{i}}}}$. If
we are unwilling or unable to specify prior covariances here, we may
obtain the $\cov{\bel}{\Ef{D}{\M{\class{i}}}}$ using a similar sampling
scheme to that described in Section \ref{algorithm}.

\section{Application: calibrating an expensive ocean model}\label{application}
We return to the ocean model calibration problem described in Section
\ref{motivating}, with the model established by Figure \ref{tBN}. We
first describe the ocean model and the data we will use to calibrate it
before returning to the problem of establishing judgements for prior
and likelihood that will enable the formal Bayesian analysis to
proceed. We will also describe a collection of alternatives and perform
the posterior belief assessment using the methodology described in
Sections \ref{sect.reified} and \ref{exchangeability}.

\subsection{The ocean model and the data}\label{data}
The NEMO ocean model \citep{madec08} simulates the global ocean for
given atmospheric forcing files. The $2^\circ$ ORCA configuration that
we use here divides the globe into approximately $2^\circ$
latitude/longitude grid boxes and has 32 depth levels. The model is
extremely computationally expensive, taking approximately 7.5 hours to
complete 30 years of model time on the UK supercomputer ARCHER. For any
setting of the parameters, we run the model for 180 model years under
climatological forcing (atmospheric forcing designed to replicate an
``average year'' in the late 20th century), and perform our analysis on
the results of the last 30 years of the simulation.
\par
\begin{figure}[ht]
\includegraphics{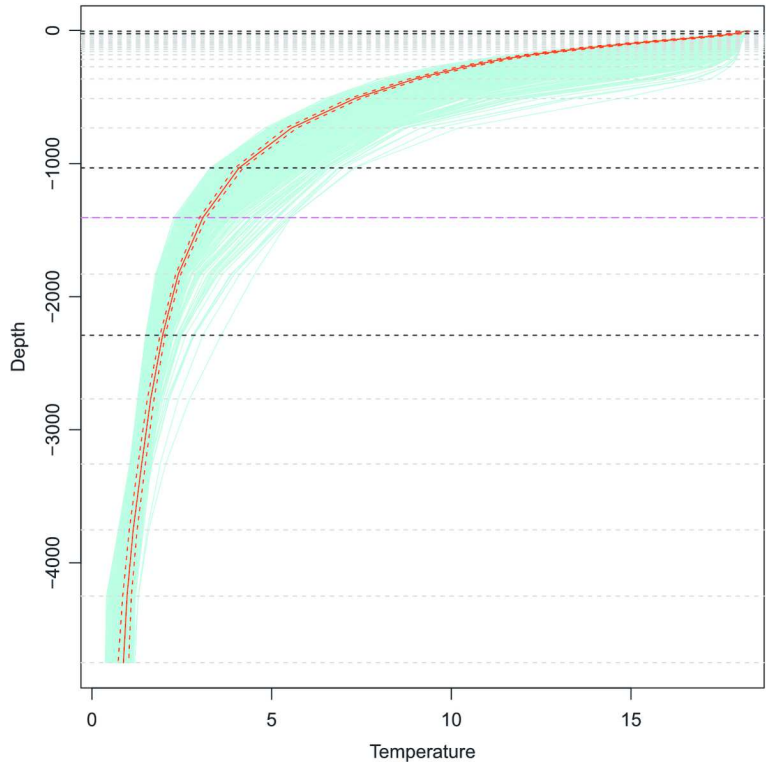}
\caption{Global mean temperatures for each of the 400 members of our
ocean model ensemble (cyan lines) plotted against model depth. The red
solid line represents observed global mean temperature and the dashed
red lines indicate $\pm2$ standard deviations of the error on the
observations. Horizontal dashed lines indicate the model depth levels,
with the black dashed lines those depths used to calibrate the model
and the magenta long-dashed line the prediction depth.}\label{ssts}
\end{figure}

The initial parameter space, containing 1 switch variable with 2
settings and 20 continuous parameters defined on a hypercube with the
ranges for each input ($a_{x^*}^{[i]}$ and $b_{x^*}^{[i]}$ in our
Bayesian network) elicited by the first author from the developer of
the code, Gurvan Madec, is explored using a 400 member $k$-extended
Latin Hypercube design (LHC), a Latin Hypercube of size 400 constructed
based on an initial LHC and extended $k-1$ times so that the design
comprises $k$ smaller LHCs of size $400/k$ (in this case $k$ was 25).
The design method and rationale for this choice applied to the $2^\circ
$ NEMO model is discussed in \cite{williamson14}. Figure \ref{ssts}
plots global mean temperature as a function of depth with cyan lines
representing the $400$ model runs and the red solid line representing
real world data \citep[using the EN3 dataset,][]{inglebyhuddleston07}.
Dashed red lines represent $\pm2$ standard deviations of the model
error calculated from \cite{inglebyhuddleston07}. The dashed horizontal
lines represent the 31 model depths.
%
\par
To illustrate a posterior belief assessment, we will calibrate the
$2^\circ$ NEMO using observations at 4 depths (illustrated by the black
dashed lines) in order to predict the temperature at a new depth for
which we pretend not to have observed data (the magenta long-dashed
line). Though clearly illustrative (the goal of this paper is not to
provide a detailed ocean reanalysis based on NEMO), this analysis will
be sufficiently complex to highlight the features of our methodology
and to provide a roadmap for posterior belief assessment in the
calibration of computer models.
\subsection{Alternative judgements}\label{nroy.sect}
Even with the detailed modelling in Section \ref{motivating}, there is
a great deal left to specify before a Bayesian analysis can be
completed. The remaining choices for prior and likelihood are part of
one member of a collection of judgements from $\jdg{0}, \jdg{1}, \ldots
$. We describe the elements of $\jdg{0}$ and list the types of possible
alternatives to each element we consider and the rationale behind
choosing them. We will then divide the types of alternatives into a
number of co-exchangeable classes, following Section \ref
{exchangeability}, and select a number of representatives for each
class. \par We begin with judgements concerning the prior and
likelihood of the emulator component, starting with the type and number
of terms in response surface component $g_d(x)$. Our preference is to
fit complex response surfaces with $g_d(x)$ containing complex
polynomial terms in $x$ chosen using some model selection algorithm.
There is considerable literature adopting this approach
\citep[\xch{cf.}{c.f.}][]{craigetal01, cumminggoldstein09, williamsonetal12, vernonetal10,
kaufmanetal11, williamsonetal14}, which is based on a philosophy that
global effects can be captured using the response surface, allowing
local effects to be captured by the Gaussian process residual. We use a
stepwise selection routine described in \cite{williamsonetal13} similar
to that used in \cite{sextonetal11} that aims to add as many variance
describing terms to each $g_d(x)$ as possible and then deletes those
that describe the least variance, one at a time, until there are fewer
terms than a maximum allowed number (specified as a percentage of the
number of degrees of freedom available) and the explained variance of
the best current best candidate for deletion is above a certain
threshold (usually $0.1\%$). Our default choice is to use $10\%$ of the
available degrees of freedom allowing up to 40 terms in each $g_d(x)$.
This is part of $\jdg{0}$.
\par
As well as the potential to choose fewer terms using our model
selection routine, we also consider other approaches that are popular
within the literature. A popular choice is $g_d(x) = (1, x_1, \ldots,
x_r)^T$ \citep{andrianakischallenor12, leeetal13}, as is the choice
made in the original computer experiments paper of $g_d(x) = 1$ \citep
{sacksetal89}. We view these two as not ruled out alternatives
reflecting, for example, a difference in expert opinion in the form of
the likelihood. In this study, as in most applications with emulators
\citep[\xch{cf.}{c.f. }][]{haylockohagan96, kennedyohagan01, santneretal03}, we
use the reference prior for response surface coefficients so that $\pi
(\beta_d, \sigma_f^2) \propto1/\sigma_f^2$.
\par
Popular choices of correlation function $R(|x-x'|;\kappa_{1:r})$ in the
literature include power exponential forms and Mat\'{e}rn forms. Power
exponential correlation functions take the form
\begin{displaymath}
R(|x-x'|;\kappa_{1:r}) = \prod_{i=1}^r\exp\{-\kappa_i|x_i-x'_i|^p\}
\end{displaymath}
with the special case $p=2$ being the Gaussian covariance function.
Though this is popular, \cite{bayarrietal07} suggested $p=1.9$ was a
better alternative as the residual under the Gaussian correlation
function was often too smooth; and others favour Mat\'{e}rn forms,
which are popular in spatial statistics. The two Mat\'{e}rn functions
we consider are chosen for computational convenience and have
\begin{displaymath}
R(|x-x'|;\kappa_{1:r}) = \prod_{i=1}^r(1 + \sqrt{3}\kappa
_i|x_i-x'_i|)\exp\{-\sqrt{3}\kappa_i|x_i-x'_i|\}
\end{displaymath}
and
\begin{displaymath}
R(|x-x'|;\kappa_{1:r}) = \prod_{i=1}^r(1 + \sqrt{5}\kappa_i|x_i-x'_i| +
\frac{1}{3}\sqrt{5}\kappa_i|x_i-x'_i|^2)\exp\{-\sqrt{5}\kappa
_i|x_i-x'_i|\}.
\end{displaymath}
In $\jdg{0}$, we select the power exponential form with $p=1.9$ as this
is the form we are most familiar with. The other forms (including
powers $p$ closer to 1) change both the prior and the likelihood and
are considered as alternatives.
\par
We used the MATCH elicitation tool \citep{morrisetal14} to select
values for the hyper-priors $a_{\kappa}^{1:r}$ and $b_{\kappa}^{1:r}$,
selecting values of $2.9$ and $5,$ respectively, for the hyper
parameters in each dimension for $\jdg{0}$. We do not rule out
alternatives with the same prior mean for each $\kappa$, but allow the
variance to be up to twice as large or as low as 5 times smaller. As
the nugget in this paper represents the proportion of the residual
variability that is due to the initial condition uncertainty in the
model, we base our prior judgements for $\nu$ on how much of the signal
in the data we believe our response surface will capture. Given
$g_d(x)$ the emulator is fitted using the joint Bayesian update
described in \cite{haylockohagan96}. However, we select $g_d(x)$ using
a stepwise selection routine (described above) based on ordinary least
squares. Following selection of $g_d(x),$ we use the $R^2$ value of the
last fit in the stepwise routine to judge how much of the residual
variability should be uncorrelated based on 7 separately elicited
scenarios. For example, the first scenario is that the OLS fit captures
more than $95\%$ of the variability in the ensemble. In this case we
believe the residual is mainly uncorrelated. We use the MATCH tool to
elicit a $\mathrm{Be}(3.8,1.7)$ distribution for $\nu$ that is skewed
towards high values. The other scenarios are $R^2>0.9$, $R^2>0.8$,
$R^2>0.7$, $R^2>0.6$, $R^2>0.4$ and $R^2\leq0.4$, with MATCH elicited
priors $\mathrm{Be}(2.3,1.7)$, $\mathrm{Be}(2,1.5)$, $\mathrm
{Be}(1.5,1.5)$, $\mathrm{Be}(1.6,2.1)$, $\mathrm{Be}(1.8, 3.8)$ and
$\mathrm{Be}(1.4, 3.1),$ respectively. These judgements are in $\jdg
{0}$. We consider alternative judgements that respect the mean of
distribution for $\nu$ but inflate the variance by multiplying both
$a_{\nu}$ and $b_{\nu}$ by a scalar as large as $4$.
\par
As we have $\sigma_e^2$ (as shown by the red dashed lines on Figure \ref
{ssts}), all that remains is to specify hyper-parameters controlling
the model discrepancy. For $\zeta$ we choose a prior for the roughness
length that reflects our belief that if the model is biased at one
depth, it is very likely to be similarly biased at nearby depths. We
specify $a_{\zeta} =1$, $b_{\zeta} = 7$ to reflect this. We do not
explore changes to these judgements, however, we concentrate on the
impact of specifying different discrepancy variance distributions. Our
current judgements for this are $a_{\eta} = 1000$, $b_{\eta}= 6.8$ $\in
\jdg{0}$. This gives an expectation for $\sigma_{\eta}^2$ equivalent to
the observation error and is a value we have considered as the
modeller's tolerance to error in our work with NEMO \citep[see][for
discussion]{williamsonetal14b}.
\par
Retaining the same mean value, we consider two alternative types of
discrepancy judgement. The first type we describe as ``medium''
discrepancy where the hyper-parameters are in the same proportion but
roughly an order of magnitude smaller so that the distribution on
$\sigma_{\eta}^2$ has a slightly larger variance and is a little more
skewed to higher values. The second is ``large'' discrepancy, which is
similar, but with 2 orders of magnitude difference in the values of
$a_{\eta}$ and $b_{\eta}$ so that $\sigma_{\eta}^2$ has a much larger
variance and is very highly skewed towards higher values. Unlike some
of the alternatives to the likelihood described earlier in this
subsection, these alternatives represent quite different perspectives
on the relevance of the NEMO ocean model to the actual ocean (for
global mean temperature at least). Hence, we might view these
alternatives as representative of different levels of confidence held
by different ocean modellers about the NEMO model.
\par
Before describing our co-exchangeable classes of alternative
judgements, it is worth making a brief comment about the way the
emulators are fitted and the calibration performed. Though we do run
MCMC to obtain posterior judgements about $y(d)$, we follow the
standard approach in calibration and first fix the correlation and
nugget parameters after conditioning on the ensemble \citep
{kennedyohagan01}. This means first fitting the emulator, then using it
in the calibration to the data. We choose the MAP estimates of $\kappa
_{1:r}$ and $\nu$ instead of maximum likelihood estimates, as these
account for our prior modelling. MAP estimates are obtained using
simulated annealing.
\subsection{Co-exchangeable classes of not ruled out alternatives}\label
{nroy.alternatives}
We will be interested in $\E{y(d_5)|z(d_{1:4})}$, the expected global
mean temperature at around 1405 metres depth given the four
observations we have taken. For many of the alternative modelling and
prior choices we have described there may well be an impact upon this
posterior expectation, however, due to the complexity of the model it
is not clear to us, a priori, how any impacts might be related.
However, we do have views regarding how the type of response surface we
fit and how the different discrepancy choices might lead us to quite
different analyses. Hence our division of the collection of possible
alternative posterior expectations under different judgements will be
based on these 2 parts of our modelling.
\par
Considering the type of response surface fitted first, we view the
alternatives as forming 3 distinct impacts upon our analysis. In the
case of complex mean functions fitted using our stepwise routine, we
consider that much of global signal will be captured by $\beta
^Tg_d(x)$, allowing the Gaussian process $\epsilon_d(x)$ to capture
only local departures from this surface. Hence posterior estimates of
$\sigma_f^2$, $\nu$ and $\kappa_{1:r}$ are likely to be much different
than in the cases where either a linear or constant response surface is
fitted. In the case of complex mean functions, as the $\epsilon_d(x)$
process will be more tuned to local deviations from the response
surface rather than the behaviour of the function in the whole model
parameter space, we believe it will be more accurate than the other
emulators in regions of parameter space that are close to our design
points. If we have any design points that are quite close to $x^*$, we
expect $\E{y(d_5)|z(d_{1:4})}$ to be captured more accurately using
this emulator than the other types.
\par
However, we have no strong views on how reducing the number of allowed
degrees of freedom to be spent on fitting the stepwise regression will
impact upon this argument. As we only allow a minimum of $5\%$ of the
degrees of freedom to be spent on the response surface (in this case
the same number of degrees of freedom spent using the linear response
surface with all of the parameters), we believe that the major global
non-linearities will be captured by $\beta^Tg_d(x)$ anyway. Therefore,
even though $\epsilon_d(x)$ may contain more of the global signal than
if a higher number of degrees of freedom were spent, it may not (some
of the additional terms retained in $\beta^Tg_d(x)$ when a higher
percentage of our degrees of freedom are spent may in fact be capturing
local signal through high order interactions between inputs).
\par
We view the cases of the linear (in the model inputs) and constant
response surfaces to have distinct impacts upon our posterior
expectation, both from the complex mean case (for the reason given
above) and from each other, as we now explain. In the case of the
emulator with linear mean, the posterior uncertainty of the emulator
will increase asymptotically as we move outside of the convex hull of
design points (as it would with the complex mean), however, the same is
not true of the constant mean emulator. Hence, we may view this
increased uncertainty outside of the convex hull of design points as
contributing to our posterior expectation in some way. Perhaps, for
example, $x^*$ is outside of that convex hull, leading to our range of
possible values of the computer model at $x^*$ to be quite different in
the two cases. Or perhaps the model outputs that are possible in the
linear mean case, but not in the constant mean case, bias the posterior
expectation in one particular direction (perhaps they allow much higher
temperatures, but not lower).
\par
So far, we have defined 3 distinct classes of judgements based only on
our choice of emulator. The choice of discrepancy variance
parametrisation (``standard'', ``medium'' or ``high'') will also impact
upon our prior judgements for $\E{y|z}$. High values of model
discrepancy variance reduce the impact of information from the ocean
model on our posterior expectations for ocean temperature. If the
discrepancy is high enough, it can effectively ``take over'' as there
is little information in the data about the ocean model parameters so
that our prediction is driven by the Gaussian process discrepancy
fitted as (approximately) a deviation from the ensemble mean. We do not
view the ``medium'' case as being different enough from the
``standard'' to change our beliefs about $\E{y|z}$. The medium case has
a wider range for $\sigma^2_{\eta}$ and is slightly more skewed towards
higher values, but only slightly and not enough to give us any strong
views. However, we view the ``high'' case, which is much more skewed and
allows much higher discrepancy variances, to be distinct from these
other two.
\par
Combinations of our two distinct classes driven by discrepancy
judgements and 3 distinct classes based on the type of emulator lead to
6 classes of potential judgements. Within each class we have all
choices of covariance function and hyper-priors for the emulator
residual. By just considering our alternative hyper prior choices
$a_{\kappa}^{1:r}$ and $b_{\kappa}^{1:r}$, (any values such that the
prior mean of the half length correlation is the same but that the
variance can be up to twice as large or as low as 5 times smaller), we
can see that there are infinitely many alternative judgements within
each class ($a_{\kappa}^{1:r}$ and $b_{\kappa}^{1:r}$ are continuous).
Within each class we have no a priori views on how the means, variances
or covariances of $\E{y|z}$ would differ with each possible
alternative. We also have no views on any specific member dependent
differences between $\E{y|z}$ calculated from different classes. Hence
we assume that the 6 defined classes $\class{1}, \ldots, \class{6}$ are
co-exchangeable.
\par
\subsection{Posterior belief assessment}
Let $\jdg{i_j}$ represent a collection of modelling judgements from
$\class{i}$, then (\ref{coex1}) gives
\begin{displaymath}
\E{y(d_5)|z(d_{1:4});\jdg{i_j}} = \M{\class{i}} + \R_{j}(\class{i}),
\end{displaymath}
and we can apply (\ref{soeM}), (\ref{soeR}), (\ref{soeCovs}), (\ref
{coexcov}) and (\ref{coexcov2}) to compute $\Ef{D}{\M{\class{i}}}$ for
$i=1, \ldots,6$, with \xch{$D=\{\mathrm{E}[y(d_5)|z(d_{1:4});\jdg{n(1)_{m(1)}}], \ldots,
\mathrm{E}[y(d_5)|z(d_{1:4});\jdg{n(N)_{m(N)}}]\}$}{$D =
\{\E{y(d_5)|z(d_{1:4});\jdg{n(1)_{m(1)}}}, \ldots,
\E{y(d_5)|z(d_{1:4});\jdg{n(N)_{m(N)}}}\}$} and with\break $\jdg
{n(1)_{m(1)}}, \ldots, \jdg{n(N)_{m(N)}}$ chosen modelling choices from
each class. Defining $\G$ as in Section \ref{exchangeability}, then in
order to perform the posterior belief assessment and compute $\Ef{\G
}{y(d_5)}$ via (\ref{reified}) we require $\E{y(d_5)}$,
$\cov{y(d_5)}{\G }$, $\var{\G}$, $\E{\G}$ and a means to computing $\G$
from $D$ which, further, requires $\E{\M{\class{i}}}$,
$\cov{\M{\class{i}}}{\M{\class {j}}}$ and $\var{\R_j(\class{i})}$ for
each $i=1, \ldots, 6$.
\par
We first select the further Bayesian analyses we will do to obtain $D$,
by choosing a representative set of judgements from each class. We do
this by first selecting a number of distinct alternative choices from
those not ruled out in Section \ref{nroy.alternatives} for each
prior/likelihood decision in our modelling. These are as follows: for
classes involving complex response surfaces, we allow both $10\%$ and
$5\%$ of the degrees of freedom to be spent in model selection. For the
``medium'' and ``high'' discrepancy settings we divide both the
standard choices of the hyper-parameters by 10 and 100, respectively.
For all classes, we allow both classes of Mat\'{e}rn, the Gaussian
covariance function and the power exponential with $p=1.5$. We multiply
all hyper-parameters of $\kappa_{1:r}$ by 5 or 0.5, respectively, and
allow the nugget hyper-parameters to be multiplied by 4. Each
combination of alternative choices and our standard choices is also allowed.
\par
Instead of running alternative Bayesian analyses for each possible
combination, we take advantage of co-exchangeability, sampling $32$
choices from $\class{1}$ (the largest class) and $8$ from each of the
remaining classes. Each Bayesian analysis builds the emulator
determined by our modelling choices, then uses a random walk metropolis
hastings algorithm with $21,000$ samples, discarding 1000 as burn in
and thinning every $20$, to obtain a sample from the posterior
distribution implied by our modelling choices from which an expectation
can be computed. This gives $D$, however, in order to compute $\G$ and
$\Ef{\G}{y(d_5)}$ we require each of the other ingredients described
above. We obtain these here by sampling.
\par
We use the sampling scheme described in Section \ref{algorithm}. The
first step of this process was to be able to form distributions $(\E
{y}, \var{y})$ and $(\E{z|y},\var{z|y})$ and to be able to sample
values of $y$ then $z$. Given $y$, we can easily obtain $z$ using
$\sigma_e^2$ and our statistical model. In an application where our
judgements would be crucial in informing scientists about global
temperatures in the real ocean, rather than one conducted for
illustrative purposes, we would prefer to attempt to elicit means and
variances for $y$ directly. However, as this example is illustrative,
we do not have the expertise available to do this. We therefore sample
values of $y$ by sampling one of our sets of judgements at random and
forming the emulator $\hat{f}(x)$, sampling an $x^*$, drawing a sample
from $\hat{f}(x)$ at this sample and adding it to a sample from the
discrepancy distribution.
\par
Having generated a $y$ and $z$ pair, we run a similar MCMC calculation
to that described above in order to generate a sample of $D$ under the
sampled $y$ and $z$ in the manner described in Section~\ref{algorithm}
step 2. To make things computationally more feasible we run shorter
simulations (6000 MCMC steps) with the same number discarded as burn
in, thinning every 10. We took 13,000 samples of the $72$ elements of
$D$ for different $y$ and $z$ pairs like this using the condor cluster
at Durham University Maths department to complete all calculations
within a day.
\par
We first use these samples to assess $\E{\M{\class{i}}}$, $\cov{\M
{\class{i}}}{\M{\class{j}}}$ and $\var{\R_j(\class{i})}$ for each $i=1,
\ldots, 6$. We obtain the first two of these three quantities for each
class by taking the mean over all experiments for each member of $D$
and by computing the mean and variances of these within each class. The
latter is assessed by taking the mean of the variances over all
experiments for each class and subtracting the assessment for the
variance of $\M{\class{i}}$. For simplicity, we assume $\cov{\M{\class
{i}}}{\M{\class{j}}}=0$ for $j\neq i$ as each had a small variance so
was treated as known. When conducting these samples, an important part
of any application of our methodology to a real world problem would be
to perform diagnostics, using the samples, on our exchangeability
judgements. Whilst we do not have the space to address this here, we
comment that the sampling method offers the chance to examine the
validity of our statistical modelling without confronting it with our
only real data.
\par
For each sample, we now compute the sample Bayes linear adjusted
expectation $\Ef{D}{\M{\class{i}}}$ using the quantities derived as
above, and use them to form sample values of $\G$ as in step 3 of our
algorithm. We now have coherent sample values of $y$ and $\G$ and use
them to establish all quantities required to compute $\Ef{\G}{y}$ via
(\ref{reified}) as described in step 4. We then compute the actual $\Ef
{D}{\M{\class{i}}}$ based on our observation driven MCMC calculations,
form $\G$ and compute $\Ef{\G}{y}$. The values of $\Ef{\G}{y(d_5)}$ and
$\E{y(d_5)|z(d_{1:4});\jdg{0}}$ are $2.921^\circ$C and $2.951^\circ$C,
respectively. The adjusted variance in $y$ following each analysis,
$\var{y(d_5)} - \var{\Ef{\G}{y(d_5)}}$ and $\var{y(d_5)} - \var{\Ef{\G
_1}{y(d_5)}}$ was $0.0262$ and $0.0226$. As described in Section \ref
{sect.reified}, part (i) of the theorem implies that a lower bound on
the amount of uncertainty reduction in $y(d_5)$ achieved performing a
posterior belief assessment (computed by taking the ratio of these two
adjusted variances) is $14.0\%$.
\par
We note here that our choice of how many samples to take from each
class was largely arbitrary in this problem. However, if we built in a
way of using elements of $D$ to learn about $\var{\R_j(\class{i})}$
before then computing $\Ef{D}{\M{\class{i}}}$ as part of what \cite
{goldsteinwooff07} called a two stage Bayes linear analysis, we would
be able to use this information to see how our uncertainty in $\M{\class
{i}}$ was reduced using $D$. The idea would be to use information from
this procedure to decide how many samples would be required from each
class. Developing a two stage procedure for posterior belief assessment
would be technically challenging and is a possible avenue of future research.
\par
The posterior belief assessment for $y$ is a constant plus a linear
combination of the elements of $\G$. Hence it is instructive to observe
what the coefficients of this linear combination are. We show these in
Table \ref{table1}. The larger coefficients for the three classes
corresponding to the high discrepancy case indicate that this judgement
may represent the main difference between our expectation following
posterior belief assessment and that conditioned only on our original
judgements. Though the goal of this analysis is not model selection,
these coefficients do offer some a posteriori information about the
relative merit of each class of model. If model selection was of
interest, any particularly influential classes might be further
explored by looking at the particular linear combination of the
elements of $D$ that influence its adjusted expectation to see if there
were any particularly influential conditional expectations. We note,
however, that the elements of $\G$ (and of $D$) are not orthogonal so
any interpretation as to the relative merits of any class of model (or
individual model) on the posterior belief assessment is not straightforward.
\par
The adjusted variance, $\var{\Ef{\G}{y}}$ is such that $\E{y|z;\jdg
{0}}$ is not even $1$ standard deviation away from $\Ef{\G}{y}$. Hence,
in this example, our modified posterior judgements are close to those
from the initial Bayesian analysis under $\jdg{0}$. This is not
strictly a robust Bayesian analysis, however, given the number of
alternative modelling choices explored, we would view our conclusions
as robust to both the prior and likelihood choices we have made. Though
there are obvious similarities with sensitivity analysis here, we note
that our inference actually combines all of the information from the
sensitivity study, rather than only using the extra information as a
sanity check. Further, if our analysis suggested, for example, that our
assessment was sensitive to the samples we had taken, that would only
suggest performing further analyses in order to reduce $\var{\Ef{\G
}{y}}$, rather than requiring us to think more carefully about certain
aspects of the prior or likelihood.
%
\begin{table}[ht]
\tabcolsep=1.6pt
\begin{tabular*}{\textwidth}{@{\extracolsep{\fill}}c|c|c|c|c|c|c@{}}
$\E{y|z;\jdg{0}}$ & $\Ef{D}{\M{\class{1}}}$ & $\Ef{D}{\M{\class{2}}}$
& $\Ef{D}{\M{\class{3}}}$ & $\Ef{D}{\M{\class{4}}}$ & $\Ef{D}{\M{\class
{5}}}$ & $\Ef{D} {\M{\class{6}}}$ \\ \hline
0.133 & 0.383 &$-$0.078 & 0.361 & 1.140 &1.024 & $-$0.842
\end{tabular*}
\caption{Coefficients from $\Ef{\G}{y}$ for each of the elements of
$\G$.}\label{table1}
\vspace*{-6pt}
\end{table}

\section{Discussion}\label{discussion}
Performing a truly subjective Bayesian analysis (where all elements of
the prior modelling are held subjective beliefs, leading to a clear
interpretation of the posterior distribution) in problems requiring
complex statistical models, would be extremely challenging, expensive
and time consuming. This would be the case even if elicitation of
high-dimensional joint distributions, even for non-standard forms, was
well supported and routine in Bayesian analysis. However, the case
where we require posterior beliefs for every, or even very many,
combinations of quantities in our statistical model is perhaps rare at
best, though likely non-existent. Posterior belief assessment is a
methodology that allows us to use the full Bayesian machinery in order
to obtain beliefs for key quantities of interest using as many
alternative forms of prior and likelihood modelling judgements as we
are prepared to consider as representative of our judgements regarding
the problem structure and our beliefs before undertaking the analysis.
\par
We showed that, as a consequence of temporal sure preference, the
posterior belief assessment is closer to the prevision we would specify
at time $t$, than a full Bayesian analysis on our initial/default set
of modelling judgements. We showed that we expect a posterior belief
assessment to resolve more of our uncertainty about those key
judgements of interest than a standard Bayesian analysis.
\par
We argued that posterior belief assessment was a powerful and tractable
alternative to traditional robust Bayesian analyses (that attempt to
proceed analytically) when of interest is how robust any key judgements
might be to plausible alternatives to modelling. Whilst a traditional
analytic robust Bayesian analysis might not be possible in most
applications with complex statistical models, we can perform a
posterior belief assessment as long as we are able to repeat the
Bayesian sampling computation (either in parallel or otherwise) for a
finite set of alternative judgements. Our approach both gives us
posterior judgements about quantities of interest, the difference
between these and a one-off Bayesian analysis, and information as to
which alternative prior modelling choices contribute to these and by
how much.
\par
Though tractable, particularly when compared to robust Bayes in complex
statistical models, a posterior belief assessment requires a number of
additional full Bayes analyses, that will likely involve
computationally expensive sampling procedures to be repeated. If it
would exhaust the analyst's computational resources to perform one of
these analyses, posterior belief assessment will not be feasible.
However, in most cases, either additional time or access to distributed
computing facilities will be available and will allow the extra
calculations to be done. In some cases, an alternative set of
judgements may lead to a statistical model for which the Bayesian
machinery has not yet been adequately developed to sample from, or that
would involve a great deal of extra effort to implement. In these
cases, it may be that such models belong to a co-exchangeable class
containing models that are easier to sample from. If so, we run the
simpler alternatives from this class and proceed as discussed in
Section \ref{exchangeability}. If difficult to implement alternatives
are not second order exchangeable with simpler alternatives, then
further methodology may be required to perform a posterior belief assessment.
\par
The principal challenge when performing a posterior belief assessment
is in considering all possible alternatives of prior and likelihood to
the original set of judgements $\jdg{0}$ that you would be unwilling to
rule out based on your current understanding. As illustrated by our
application to the ocean model, in complex models the ways in which the
model can be plausibly (in your view) changed can quickly grow and be
difficult to think about and to then group into co-exchangeable
classes. There are many statistical models being developed that are far
more complex than our ocean model example (but that can be readily
sampled from using MCMC), so this task could seem far more daunting
than in our case in some applications. However, in subjective Bayes
there is no such thing as a free lunch and, making meaningful belief
statements updated properly by available data, necessarily involves
careful thought so that their meaning is well understood following the analysis.
\par
The other main challenge in performing a posterior belief assessment
involves obtaining the quantities required to compute $\Ef{\G}{\bel}$,
namely $\E{\bel}$, $\cov{\bel}{\G}$, $\var{\G}$ and $\E{\G}$. If these
could not be elicited directly using partial prior specification
methods discussed in chapter 2 of \cite{goldsteinwooff07}, a
combination of expert judgement (for quantities only involving $\bel$),
and a sampling scheme we outlined in Section \ref{algorithm} could be
used. The sampling method required many many more Bayesian computations
to be performed (perhaps in more approximate form, e.g. MCMC with fewer
samples), and may be infeasible if access to distributed computing is
unavailable or in particularly complex problems. In particular, the
number of Bayesian calculations required for a posterior belief
assessment using our sampling algorithm will require a certain amount
of automation in order to be feasible. Distributed computing programs
such as condor make this feasible, as highlighted in our application.
Our application was reasonably complex, yet we managed to perform
$13000 \times72$ separate Bayesian analyses based on MCMC using the
automation provided by the condor program and a $1000$ core cluster at
Durham University in only 1 day. Developing tools to allow elicitation
of these quantities directly, or developing alternative methods of
deriving them through further, less computationally burdensome,
calculations could be an important avenue of further investigation in
this area.



%

%
\begin{acknowledgement}
Danny Williamson is funded by an EPSRC fellowship grant number
EP/K019112/1. We wish to thank an associate editor and two referees for
their thoughtful and detailed comments regarding the paper. Their
suggestions have helped to improve it a great deal. We also thank Adam
Blaker and Bablu Sinha for running the NEMO ocean model to generate the
ensemble that appears in our application.
\end{acknowledgement}


\begin{thebibliography}{}

\bibitem[Andrianakis and Challenor(2012)]{andrianakischallenor12}
Andiranakis, I. and Challenor, P.~G. (2012), The effect of the nugget on
Gaussian process emulators of computer models. \emph{Computational
Statistics and Data Analysis}, 56: 4215--4228.
\bid{doi={10.1016/j.csda.2012.04.020}, issn={0167-9473}, mr={2957866}}
\bptok{addids}%
\endbibitem

\bibitem[Bayarri et~al.(2007)]{bayarrietal07}
Bayarri, M.~J., Berger, J.~O., Cafeo, J., Garcia-Donato, G., Liu, F., Palomo,
J., Parthasarathy, R.~J., Paulo, R., Sacks, J., and Walsh, D. (2007),
Computer model validation with functional output, \emph{The
Annals of Statistics}, 35, 1874--1906.
\bid{doi={10.1214/009053607000000163}, issn={0090-5364}, mr={2363956}}
\bptok{addids}%
\endbibitem

\bibitem[Berger(1994)]{berger94}
Berger, J.~O. (1994), An overview of robust Bayesian analysis,
\emph{Test}, 3(1), 5--124.
\bid{doi={10.1007/BF02562676}, issn={1133-0686}, mr={1293110}}
\bptok{addids}%
\endbibitem

\bibitem[Berger(2006)]{berger06}
Berger, J.~O. (2006), The case for objective Bayesian analysis.
\emph{Bayesian Analysis}, 1(3), 385--402.
\bid{issn={1936-0975}, mr={2221271}}
\bptok{addids}%
\endbibitem

\bibitem[Craig et al.(2001)]{craigetal01}
Craig, P.~S., Goldstein, M., \xch{Rougier, J.~C.}{C., R.~J.}, and Seheult, A.~H. (2001),
Bayesian forecasting for complex systems using computer simulators,
\emph{Journal of the American Statistical Association}, 96, 717--729.
\bid{doi={10.1198/\\016214501753168370}, issn={0162-1459}, mr={1946437}}
\bptok{addids}%
\endbibitem

\bibitem[Cumming and Goldstein(2009)]{cumminggoldstein09}
Cumming, J.~A. and Goldstein, M. (2009), Small sample designs for
complex high-dimensional models based on fast approximations,
\emph{Technometrics}, 51, 377--388.
\bid{doi={10.1198/TECH.2009.08015}, issn={0040-1706}, mr={2756474}}
\bptok{addids}%
\endbibitem

\bibitem[de Finetti(1974)]{definetti74}
de Finetti, B. (1974), \emph{Theory of Probability, Volume 1},
Wiley, New York.
\bptok{addids}%
\endbibitem

\bibitem[de Finetti(1975)]{definetti75}
de Finetti, B. (1975), \emph{Theory of Probability, Volume 2},
Wiley, New York.
\bptok{addids}%
\endbibitem

\bibitem[Gelman et al.(2004)]{gelmanetal04}
Gelman, A., Carlin, J.~B., Stern, H.~S., and Rubin, D.~B. (2004), \textit
{Bayesian Data Analysis, second edition}, Chapman and Hall/CRC Texts in
Statistical Science.
\bid{mr={2027492}}
\bptok{addids}%
\endbibitem

\bibitem[Goldstein(1997)]{goldstein97}
Goldstein, M. (1997). Prior inferences for posterior judgements. In:
\emph{Structures and Norms in Science}, M.~C.~D. Chiara et al., eds.,
Kluwer, 55--71.
\bid{mr={1819330}}
\bptok{addids}%
\endbibitem

\bibitem[Goldstein(2006)]{goldstein06}
Goldstein, M. (2006), Subjective Bayesian analysis: Principles and
practice. \emph{Bayesian Analysis}, 1(3), 403--420.
\bid{mr={2221272}}
\bptok{addids}%
\endbibitem

\bibitem[Goldstein(2011)]{goldstein11}
Goldstein, M. (2011). External Bayesian analysis for computer
simulators. In: \emph{Bayesian Statistics 9}. Bernado, J.~M. et al.,
eds., Oxford University Press.
\bid{doi={10.1093/acprof:oso/9780199694587.003.0007}, mr={3204007}}
\bptok{addids}%
\endbibitem

\bibitem[Goldstein(2012)]{goldstein12}
Goldstein, M. (2012). Observables and models: exchangeability and the
inductive argument. In: \emph{Bayesian Theory and Applications},
Damien, P. et al., eds., Clarendon Press, Oxford, 3--18.
\bid{doi={10.1093/acprof:oso/\\9780199695607.003.0001}, mr={3221155}}
\bptok{addids}%
\endbibitem

\bibitem[Goldstein and Wooff(2007)]{goldsteinwooff07}
Goldstein, M. and Wooff, D. (2007), \textit{Bayes Linear Statistics
Theory and Methods}, John Wiley and Sons Ltd.
\bid{doi={10.1002/\\9780470065662}, mr={2335584}}
\bptok{addids}%
\endbibitem

\bibitem[Haylock and O'Hagan(1996)]{haylockohagan96}
Haylock, R. and O'Hagan, A. (1996), \enquote{On inference for outputs of
computationally expensive algorithms with uncertainty on the inputs.} In:
\textit{Bayesian Statistics 5}, Bernado, J.~M., Berger, J.~O., Dawid,
A.~P., and Smith, A. F.~M., eds., Oxford University Press, \xch{629--637}{pp. 629--637}.
\bid{mr={1425432}}
\bptok{addids}%
\endbibitem

\bibitem[Higdon et al.(2008)]{higdonetal08}
Higdon, D., Nakhleh, C., Gattiker, J., and Williams, B. (2008), \enquote{A
Bayesian calibration approach to the thermal problem}, \textit{Computer
Methods in Applied Mechanics and Engineering,} 197, 2431--2441.
\bptok{addids}%
\endbibitem

\bibitem[Hoeting et al.(1999)]{hoetingetal99}
Hoeting, J.~A., Madigan, D., Raftery, A.~E., and Volinsky, V.~T. (1999),
Bayesian model averaging: a tutorial, \emph{Statistical Science},
14(4), 382--401.
\bid{doi={10.1214/ss/1009212519}, issn={0883-4237}, mr={1765176}}
\bptok{addids}%
\endbibitem

\bibitem[Ingleby and Huddleston(2007)]{inglebyhuddleston07}
Ingleby, B. and Huddleston, M. (2007), Quality control of ocean
temperature and salinity profiles -- Historical and real time data,
\emph{Journal of Marine Systems}, 65 158--175.
\bptok{addids}%
\endbibitem

\bibitem[Kaufman et al.(2011)]{kaufmanetal11}
Kaufman, C.~G., Bingham, D., Habib, S., Heitmann, K., and Frieman,
J.~A. (2011), Efficient emulators of computer experiments using
compactly supported correlation functions, with an application to
cosmology, \emph{The Annals of Applied Statistics,} 5(4), 2470--2492.
\bid{doi={10.1214/11-AOAS489}, issn={1932-6157}, mr={2907123}}
\bptok{addids}%
\endbibitem

\bibitem[Kennedy and O'Hagan(2001)]{kennedyohagan01}
Kennedy, M.~C. and O'Hagan, A. (2001), Bayesian calibration of computer
models. \emph{Journal of the Royal Statistical Society: Series B
(Statistical Methodology),} 63: 425--464.
\bid{doi={10.1111/1467-9868.00294}, issn={1369-7412}, mr={1858398}}
\bptok{addids}%
\endbibitem

\bibitem[Lee et al.(2013)]{leeetal13}
Lee, L., Pringle, K., Reddington, C., Mann, G., Stier, P., Spracklen,
D., Pierce, J., and Carslaw, K. (2013), The magnitude and causes of
uncertainty in global model simulations of cloud condensation nuclei,
\emph{Atmospheric Chemistry and Physics}, 13.17, 8879--8914.
\bptok{addids}%
\endbibitem

\bibitem[Liang et al.(2010)]{liangetal10}
Liang, F., Liu, C., and Carroll, R.~J. (2010), \textit{Advanced Markov
Chain Monte Carlo Methods}, John Wiley and Sons Ltd., Chichester, UK.
\bid{doi={10.1002/9780470669723}, mr={2828488}}
\bptok{addids}%
\endbibitem

\bibitem[Lindley(2000)]{lindley00}
Lindley, D.~V. (2000), The philosophy of statistics. \emph{Journal of
the Royal Statistical Society: Series D (The Statistician),} 49(3), 293--337.
\bptok{addids}%
\endbibitem

\bibitem[Madec(2008)]{madec08}
Madec, G. (2008), NEMO ocean engine. \emph{Note du Pole de Mod\'
elisation}. Institut Pierre-Simon Laplace (IPSL), France, No. 27, ISSN
No. 1288-1619.
\bptok{addids}%
\endbibitem

\bibitem[Morris et al.(2014)]{morrisetal14}
Morris, D.~E., Oakley, J.~E., and Crowe, J.~A. (2014), A web-based tool for
eliciting probability distributions from experts. \emph{Environmental
Modelling and Software}, 52: 1--4, ISSN 1364-8152, \url
{http://dx.doi.org/10.1016/j.envsoft.2013.10.010.}
\bptok{addids}%
\endbibitem

\bibitem[Oakley and O'Hagan(2010)]{oakleyohagan10}
Oakley, J.~E. and O'Hagan, A. (2010). \textit{SHELF: the Sheffield Elicitation
Framework (Version 2.0)},
School of Mathematics and Statistics, University of Sheffield (2010),
\url{www.tonyohagan.co.uk/shelf}
\bptok{addids}%
\endbibitem

\bibitem[Sacks et al.(1989)]{sacksetal89}
Sacks, J., Welch, W.~J., Mitchell, T.~J., and Wynn, H.~P. (1989),
Design and analysis of computer experiments, \emph{Statistical
Science}, 4, 409--435.
\bid{issn={0883-4237}, mr={1041765}}
\bptok{addids}%
\endbibitem

\bibitem[Santner et al.(2003)]{santneretal03}
Santner, T.~J., Williams, B.~J., and Notz, W.~I. (2003), \emph{The
Design and
Analysis of Computer Experiments}, Springer-Verlag New York.
\bid{doi={10.1007/978-1-4757-3799-8}, mr={2160708}}
\bptok{addids}%
\endbibitem

\bibitem[Savage(1977)]{savage77}
Savage, L.~J. (1977) The shifting foundations of statistics. In: \textit{Logic,
Laws and Life: Some Philosophical Complications}, R.G. Colodny, ed.,
Pittsburgh University Press, Pittsburg, 3--18.
\bptok{addids}%
\endbibitem

\bibitem[Sexton et al.(2011)]{sextonetal11}
Sexton, D.~M.~H., Murphy, J.~M., and Collins, M. (2011), Multivariate
probabilistic projections using imperfect climate models part 1:
outline of methodology, \emph{Climate Dynamics}, \bid{doi={10.1007/s00382-011-1208-9}}
\bptok{addids}%
\endbibitem

\bibitem[Vernon et al.(2010)]{vernonetal10}
Vernon, I., Goldstein, M., and Bower, R.~G. (2010), \enquote{Galaxy formation:
a Bayesian uncertainty analysis,} \textit{Bayesian Analysis,} 5(4),
619--846, with discussion.
\bid{doi={10.1214/10-BA524}, issn={1936-0975}, mr={2740148}}
\bptok{addids}%
\endbibitem

\bibitem[Williamson and Blaker(2014)]{williamsonblaker14}
Williamson, D. and Blaker, A.~T. (2014), Evolving Bayesian emulators for
structurally chaotic time series with application to large climate
models. \emph{SIAM/ASA Journal on Uncertainty Quantification}, 2(1) 1--28.
\bid{doi={10.1137/120900915}, issn={2166-2525}, mr={3283898}}
\bptok{addids}%
\endbibitem

\bibitem[Williamson et al.(2012)]{williamsonetal12}
Williamson, D., Goldstein, M., and Blaker, A. (2012), \enquote{Fast
linked analyses for scenario based hierarchies,} \textit{Journal of the
Royal Statistical Society: Series C (Applied Statistics)}, 61(5), 665--692.
\bid{doi={10.1111/\\j.1467-9876.2012.01042.x}, issn={0035-9254}, mr={2993504}}
\bptok{addids}%
\endbibitem

\bibitem[Williamson et al.(2013)]{williamsonetal13}
Williamson, D., Goldstein, M., Allison, L., Blaker, A., Challenor, P.
Jackson, L., and Yamazaki, K. (2013), History matching for exploring and
reducing climate model parameter space using observations and a large
perturbed physics ensemble. \emph{Climate Dynamics} 41:
1703--1729. \bid{doi={10.1007/s00382-013-1896-4}}
\bptok{addids}%
\endbibitem

\bibitem[Williamson(2014)]{williamson14}
Williamson, D. (2015). Exploratory ensemble designs for environmental
models using $k$-extended Latin hypercubes \emph{Environmetrics}, 26(4)
268--283.
\bptok{addids}%
\endbibitem

\bibitem[Williamson et al.(2014)]{williamsonetal14}
Williamson, D., Blaker, A.~T., Hampton, C., and Salter, J. (2014).
Identifying and removing structural biases in climate models with
history matching. \emph{Climate Dynamics}, Online First,
\bid{doi={10.1007/s00382-014-2378-z}}
\bptok{addids}%
\endbibitem

\bibitem[Williamson et al.(2015)]{williamsonetal14b}
Williamson, D., Blaker, A.~T., and Sinha, B. (2015) Statistical ocean
model tuning and parametric uncertainty quantification with NEMO, \emph
{In preparation}.
\bptok{addids}%
\endbibitem

\end{thebibliography}
\end{document}